\documentclass[11pt]{article}
\usepackage{rotcapt}

\usepackage{graphicx}
\usepackage{amsmath}
\usepackage{amsfonts}
\usepackage{amssymb}
\usepackage{graphicx}

\renewcommand{\arraystretch}{1.3}

\catcode`\@=11
\def\marginnote#1{}

\newcount\hour
\newcount\minute
\newtoks\amorpm
\hour=\time\divide\hour by60
\minute=\time{\multiply\hour by60 \global\advance\minute by-\hour}
\edef\standardtime{{\ifnum\hour<12 \global\amorpm={am}%
        \else\global\amorpm={pm}\advance\hour by-12 \fi
        \ifnum\hour=0 \hour=12 \fi
        \number\hour:\ifnum\minute<10 0\fi\number\minute\the\amorpm}}
\edef\militarytime{\number\hour:\ifnum\minute<10 0\fi\number\minute}

\def\draftlabel#1{{\@bsphack\if@filesw {\let\thepage\relax
      \xdef\@gtempa{\write\@auxout{\string
          \newlabel{#1}{{\@currentlabel}{\thepage}}}}}\@gtempa \if@nobreak
    \ifvmode\nobreak\fi\fi\fi\@esphack} \gdef\@eqnlabel{#1}}
    \def\@eqnlabel{}
\def\@vacuum{}
\def\draftmarginnote#1{\marginpar{\raggedright\scriptsize\tt#1}}

\def\draft{
  \oddsidemargin -.5truein
  \def\@oddfoot{\footnotesize \sl preliminary draft \hfil
    \rm\thepage\hfil\sl\today\quad\militarytime}
  \let\@evenfoot\@oddfoot \overfullrule 3pt
    \let\label=\draftlabel
    \let\marginnote=\draftmarginnote
  \def\@eqnnum{(\theequation)\rlap{\kern\marginparsep\tt\@eqnlabel}%
    \global\let\@eqnlabel\@vacuum}

  }

\makeatletter
\newdimen\normalarrayskip              
\newdimen\minarrayskip                 
\normalarrayskip\baselineskip
\minarrayskip\jot
\newif\ifold             \oldtrue            \def\new{\oldfalse}
\def\arraymode{\ifold\relax\else\displaystyle\fi} 
\def\eqnumphantom{\phantom{(\theequation)}}     
\def\@arrayskip{\ifold\baselineskip\z@\lineskip\z@
     \else
     \baselineskip\minarrayskip\lineskip2\minarrayskip\fi}
\def\@arrayclassz{\ifcase \@lastchclass \@acolampacol \or
\@ampacol \or \or \or \@addamp \or
   \@acolampacol \or \@firstampfalse \@acol \fi
\edef\@preamble{\@preamble
  \ifcase \@chnum
     \hfil$\relax\arraymode\@sharp$\hfil
     \or $\relax\arraymode\@sharp$\hfil
     \or \hfil$\relax\arraymode\@sharp$\fi}}
\def\@array[#1]#2{\setbox\@arstrutbox=\hbox{\vrule
     height\arraystretch \ht\strutbox
     depth\arraystretch \dp\strutbox
     width\z@}\@mkpream{#2}\edef\@preamble{\halign
\noexpand\@halignto
\bgroup \tabskip\z@ \@arstrut \@preamble \tabskip\z@ \cr}%
\let\@startpbox\@@startpbox \let\@endpbox\@@endpbox
  \if #1t\vtop \else \if#1b\vbox \else \vcenter \fi\fi
  \bgroup \let\par\relax
  \let\@sharp##\let\protect\relax
  \@arrayskip\@preamble}
\def\eqnarray{\stepcounter{equation}%
              \let\@currentlabel=\theequation
              \global\@eqnswtrue
              \global\@eqcnt\z@
              \tabskip\@centering
              \let\\=\@eqncr

 \halign to \displaywidth\bgroup
    \eqnumphantom\@eqnsel\hskip\@centering
    $\displaystyle \tabskip\z@ {##}$%
    \global\@eqcnt\@ne \hskip 2\arraycolsep
         $\displaystyle\arraymode{##}$\hfil
    \global\@eqcnt\tw@ \hskip 2\arraycolsep
         $\displaystyle\tabskip\z@{##}$\hfil
         \tabskip\@centering
    &{##}\tabskip\z@\cr}
\newfont{\hr}{msbm10}
\newfont{\ams}{msam10}

\hoffset= - 0.9in         

\def\beq{\begin{equation}}
\def\eeq{\end{equation}}
\def\ba{\beq\new\begin{array}{c}}
\def\ea{\end{array}\eeq}

\def\be{\ba}
\def\ee{\ea}

\def\N2{${\cal N}=2$}

\def\1N{${\cal N}=1$}
\def\4N{${\cal N}=4$}
\def\nn{\nonumber}

\newdimen\linethick  \linethick=0.4pt
\newdimen\hboxitspace    \hboxitspace=5pt
\newdimen\vboxitspace    \vboxitspace=5pt

\def\fr#1{%
\beq\new
\vcenter{
\hrule height\linethick
          \hbox{\vrule width\linethick
                \kern\hboxitspace
                \vbox{\kern\vboxitspace
                      \hbox{$\begin{array}{c}\displaystyle#1
         \end{array}$}%
                      \kern\vboxitspace}%
                \kern\hboxitspace
                \vrule width\linethick}%
          \hrule height\linethick}%
\eeq}

\renewcommand{\tt}[1][mer]{\hbox{\tiny{#1}}}

\newcommand{\Tr}{\mathop{\rm Tr}\nolimits}

\def\tr{{\rm tr}\,}
\def\Tr{{\rm Tr}\,}

\def\l[{\phantom.[}
\def\ldb{\left(\!\!\left(}
\def\rdb{\right)\!\!\right)}

\def\K{h}
\def\LRC{Littlewood-Richardson coefficients}
\def\CGC{Racah 
coefficients}

\hoffset= - 1.0in         

\newdimen\linethick  \linethick=0.4pt
\newdimen\hboxitspace    \hboxitspace=5pt
\newdimen\vboxitspace    \vboxitspace=5pt

\def\fr#1{%
\beq\new
\vcenter{
\hrule height\linethick
          \hbox{\vrule width\linethick
                \kern\hboxitspace
                \vbox{\kern\vboxitspace
                      \hbox{$\begin{array}{c}\displaystyle#1
         \end{array}$}%
                      \kern\vboxitspace}%
                \kern\hboxitspace
                \vrule width\linethick}%
          \hrule height\linethick}%
\eeq}

\title{{\bf Character expansion
for HOMFLY polynomials. II. Fundamental representation. Up to five
strands in braid  } \vspace{.2cm}}
\author{{\bf A.Mironov}\footnote{ {\small {\it
Lebedev Physics Institute} and {\it ITEP, Moscow, Russia}};
mironov@itep.ru; mironov@lpi.ru}, {\bf A.Morozov}\thanks{{\small
{\it ITEP, Moscow, Russia}}; morozov@itep.ru}, {\bf
And.Morozov}\thanks{{\small {\it Moscow State University} and
{\it ITEP Moscow, Russia}};
Andrey.Morozov@itep.ru}\date{ }}

\begin{document}
 \maketitle

\vspace{-6.5cm}

\begin{center}
\hfill FIAN/TD-18/11\\
\hfill ITEP/TH-50/11\\
\end{center}

\vspace{5cm}

\centerline{ABSTRACT}

\bigskip

{\footnotesize Character expansion is introduced and explicitly
constructed for the (non-colored) HOMFLY polynomials of the simplest
knots. Expansion coefficients are not the knot invariants and can
depend on the choice of the braid realization. However, the method
provides the simplest systematic way to construct HOMFLY
polynomials directly in terms of the variable $A=q^N$: a much better
way than the standard approach making use of the skein relations.
Moreover, representation theory of the simplest quantum group
$SU_q(2)$ is sufficient to get the answers for all braids with
$m< 5$ strands. Most important we reveal a hidden hierarchical
structure of expansion coefficients, what allows one to express all
of them through extremely simple elementary constituents.
Generalizations to arbitrary knots and arbitrary representations is
straightforward. }

\bigskip

\section{Introduction}

Character expansions play an increasing role in modern studies of
exact (non-perturbative) partition functions in various quantum
field theory models, from QCD \cite{charexp} to conformal field theories \cite{cft}
and matrix models \cite{mamocharexp}. They help to reveal explicit and, especially,
hidden symmetries, in particular, are important for the study of
hidden integrability properties \cite{intcharexp}, which reflect the existence of
non-linear relations between the correlation functions.
Not surprisingly, these expansions provide also a powerful tool for
explicit calculations.

So far, most applications of character expansions arise in
Yang-Mills theories, especially on the lattice \cite{latYM}, and in the theory of
matrix models \cite{mamoYM}. In \cite{AMMceI} we suggested to extend these
considerations to Chern-Simons theory of knots \cite{CS}, namely to the theory of
HOMFLY polynomials \cite{HOMFLY} and superpolynomials \cite{GSV}, the key players in the
theory of knot invariants. As explained in \cite{AMMceI}, in the
well studied case of the torus knots, character expansions are
indeed useful to explicitly express the linear and non-linear
relations between the HOMFLY polynomials in the form of the
"$A$-polynomial" difference equations and the Hirota/Plucker relations
respectively. This adds to the previous demonstration in
\cite{DMMSS} of how the character decompositions can be used to
construct generic superpolynomials and superseries for toric knots
and links, which is by itself a highly non-trivial problem. Though
all these results do not extend literally to arbitrary knots, they
provide a serious motivation for the study of character expansions
in knot theory. Perhaps surprisingly, this direction did not attract
as much attention so far as it clearly deserves. The lucky
exceptions are by-now classical papers \cite{knotcharexp}, but they do not
go too deep into the structure of the expansions and thus do not
reveal clearly its very interesting properties.

Perhaps, the reason for an insufficient attention to character
expansions in the case of knots, is that they are {\it not} knot
invariants: the expansion depends on the braid realization of the knot,
and knot equivalent braid realizations provide different expansions
of one and the same HOMFLY polynomial. However, we shall see that
even for knot invariants this technique is extremely useful: it
provides very simple formulas directly for the HOMFLY polynomials (i.e.
directly in terms of $A=q^N$ rather than for particular values of
$N$), moreover, for entire infinite series of knots (of which the
torus series is looking just a non-specific example), what is
hardly achievable for the alternative approach based on
(computer) application of the skein relations. In fact, there are
much more applications of the character
expansion \cite{AMMceI}: from associating integrable structures with
knots to effective dealing with A-polynomials \cite{Apoly}.

\bigskip

In the present paper we study the character decomposition of the HOMFLY
polynomials \cite{AMMceI},\footnote{As usual, $\sum_{Q \vdash K}$ means a sum over
all Young diagrams $Q$ of the size (the number of boxes) equal to $K$.
In what follows, we parameterize Young diagrams by a partition $Q=\{q_1\ge
q_2\ge\ldots\ge 0\}$.}
\be\label{HOexpan}
H_R^{\cal K} = \sum_{Q \vdash m|R|}
h_R^Q S_R^*
\ee
where $S^*_R$ are the Schur functions(characters of the linear groups $GL(N)$
taken at the special point of the time-variable space $p_k=p_k^*$,
see eq.(\ref{p*slice}) below,
and find expressions for the coefficients $h_R^Q$ in
the Turaev-Reshetikhin formalism \cite{TR}, i.e. in terms of the products of
quantum ${\cal R}$-matrices along the $m$-strand braid, arising in a
$2d$ projection of the knot ${\cal K}$ (from the point of view of
Chern-Simons theory this corresponds to evaluating the functional
integral in the temporal gauge ${\cal A}^0=0$, see \cite{MorSmi} for
the current status of such an interpretation). As already mentioned,
such a decomposition, and particular coefficients $h_R^Q$ depend on
the braid realization: it is enough to mention that even the number
of strands $m$ is not by itself a knot invariant. We shall see,
however, that the dependence on the choice of the braid realization
is not as strong as it could be, perhaps even some knot {\it
co}variance of the coefficient sets $\{h_R^Q\}$ can be found in this
context to substitute/generalize the knot {\it in}variance of the
HOMFLY polynomials $H_R^{\cal K}$.

The next step is to switch from the ordinary to Tanaka-Krein
representation of the ${\mathfrak{R}}$-matrices
(for a categorical approach to the construction described here see
\cite{MF}). Namely, with the $m$-strand
braid one naturally associates a decomposition of the representation
product,
\be\label{Rmdeco}
R^{\otimes m} = \oplus_{Q \vdash m|R|} {\cal
M}_Q\otimes Q
\ee
and the crucial property of ${\mathfrak{R}}$-matrices is
that they act as units in all the constituent representations $Q$.
Thus, the ${\mathfrak{R}}$-matrices are naturally projected to the "space of
representations"
\be
{\cal M} = \oplus_{Q \vdash m|R|} {\cal M}_Q
\ee
and have there a block-diagonal form
\be
\widehat{\mathfrak{R}} =
{\rm diag}\{\widehat {\mathfrak{R}}^Q\}
\ee
The matrices $\widehat{\mathfrak{
R}}^Q$ can also be diagonalized, but now one should recall that there
are actually $m-1$ different ${\mathfrak{R}}$-matrices in the braid
realization, with ${\mathfrak{R}}_\mu$ acting at the intersection of the
adjacent strands $\mu$ and $\mu+1$. So, each $\widehat{\mathfrak{
R}}^Q_\mu$ can be diagonalized, but not for all values of $\mu \in
\{1,\ldots,m-1\}$ at once. Instead, $\widehat{\mathfrak{ R}}_\mu$-matrices
for different $\mu$ are related by the conjugation,
\be
\widehat{\mathfrak{
R}}_\mu = {\widehat {\cal U}_{\mu\nu}} \widehat{\mathfrak{R}}_\nu{\widehat {\cal
U}_{\mu\nu}}^{-1}
\ee
where the "mixing" matrices ${\widehat {\cal U}}$ can be chosen orthogonal
and also have a block diagonal form
\be
\widehat{\cal U}_{\mu\nu} =
{\rm diag}\{\widehat {\cal U}^Q_{\mu\nu}\}
\ee

However, this is not the end of the story. Mixing matrices can be
further decomposed into elementary constituents, which appear to
exhibit additional universality properties and are directly related to
the Racah coefficients. It looks plausible that,
after some work, the coefficients $h_R^Q$ for generic knots can be
all expressed in an absolutely explicit form. In this particular
paper we demonstrate how all this works in the simplest, still
non-trivial case of the fundamental representation $R=[1]$ and the
small number of strands $m\leq 5$. This is actually enough to
explicitly express  {\it all} the HOMFLY polynomials for all the
knots $P_i$ from the Rolfsen table at \cite{katlas} with $P\leq 8$
through just two discrete functions $\varkappa_Q$ and $C_k =
1/[k]_q$. In
further papers of the series this result will be extended to other
representations $R$ (to the colored HOMFLY polynomials) and to the
broader braids with $m\ge 5$.

Now we explain this general scheme more concretely. First of all,
as we already mentioned, the character expansion (\ref{HOexpan}) is taken
not an {\it arbitrary} point in the space of time-variables:
it is constrained to just a $2$-dimensional slice which encodes all the
dependence on the group in the $A$-dependence,
\be
p^*_k = \frac{A^k-A^{-k}}{q-q^{-1}} = \frac{\{A^k\}}{\{q\}}
\label{p*slice}
\ee
Hereafter, we introduced a useful notation $\{x\} = x - x^{-1}$ to simplify
the formulas.
For $A=q^N$ these $p^*_k = [N]_q$, where the $q$-number is defined
$[k]_q\equiv (q^k-q^{-k})/(q-q^{-1})$.

The manifest expressions for the Schur functions $S_Q\{p^*\}$ in these special points (\ref{p*slice})
 are quite simple and generalize the standard hook formula \cite{hook}:
\be\label{he}
S_Q\{p^*\}= \prod_{(i,j)\in Q} \frac{\{Aq^{i-j}\}}{\{q^{h_{i,j}}\}} \ \ \ \
\stackrel{A=q^N}{\longrightarrow}\ \ \ \ { \prod_{(i,j)\in Q} \frac{[N+i-j]_q}{[h_{i,j}]_q} }
\ee
where $h_{i,j}$ is the hook length.

Now, there are two important facts:

\begin{itemize}
\item[{\bf (i)}]
The general expression for the HOMFLY polynomial within the
Reshetikhin-Turaev approach is given
\be
H_R^{{\cal K}}\{p^*\} = \Tr_{R^{\otimes m}} {\cal B}^{\cal K}
= \Tr_{R^{\otimes m}} \prod_s \mathfrak{R}_{\mu(s)}^\pm
\label{HBR}
\ee
as a weighted $A$-dependent trace of an element ${\cal B}^{\cal K}$
of an $m$-strand braid group, which is a product of quantum ($q$-dependent)
$\mathfrak{R}$-matrices appearing in the braid in a certain sequence
labeled by the index $s$ in (\ref{HBR}). Trace here is a weighted trace,
see eq.(\ref{trace}) below.

\item[{\bf (ii)}]
As we discussed above, the quantum $\mathfrak{R}$-matrix
acts as a $c$-number
in irreducible representations $Q$ in the decomposition (\ref{Rmdeco}).
This statement, however, requires a more careful formulation.
In (\ref{Rmdeco}) we denote by $Q$ the highest weight representation
(in practice, it is labeled by the Young diagram) so that the sum runs over highest weights, or Young
diagrams.
${\cal M}_Q$ is actually an intertwining operator.
When the representation $Q$ with the same highest weight appears several times
in the expansion of
$R^{\otimes m}$, the space of intertwining operators has a non-unit
dimension $N_{R^m}^Q = {\rm dim}\, {\cal M}_{R^m}^Q$
known as Littlewood-Richardson coefficient.
In the present paper $R$ is actually $R=[1]$,
and in what follows we often omit the subscript $R^m$.

\end{itemize}

These two well-known facts immediately lead to decomposition (\ref{HOexpan}),
with coefficients $\K_R^Q$ expressed through the eigenvalues of
$\mathfrak{R}$-matrix and the Racah coefficients,
which seem to be quite comprehensible.
Moreover, the only source of $A$ (and thus $N$) dependence is the traces
\be\label{trace}
\Tr_Q\ I \equiv \ \hbox{ordinary trace}_Q (q^\rho)^{\otimes m} = S_Q^* = S_Q\{p^*\}
\ee
where we manifestly included the factor $q^\rho$ into the definition of trace.
Therefore, the coefficients $\K_R^Q$ can be calculated for the smallest
possible group $SU(N=l(Q))$, where $l(Q)$ is the number of columns
in the Young diagram describing the highest weight of $Q$. From now on, we denote through
$Q$ both the representation and the corresponding Young diagram, hopefully this would not
cause any misinterpretation.

\bigskip

Our main goal in this text is a review of the statements (i) and (ii)
and their not-so-trivial relation to the usual {\it straightforward}
approaches to evaluation of the HOMFLY polynomials,
which do {\it not} reveal the hidden structure (\ref{HOexpan}).
The key point will be reformulation of the braid traces:
from matrices in the representation space $R^{\otimes m}$ we switch
to those in the space
${\cal M} = \oplus_Q {\cal M}^Q$
of the intertwining operators, appearing in the decomposition
(\ref{Rmdeco}).
The true meaning of (ii) is that each ${\cal M}_Q$ is preserved
by action of the $\mathfrak{R}$-matrices
and these latter can be converted into the block-diagonal matrices which
act on ${\cal M}$.
In what follows we denote $\mathfrak{R}$-matrices acting on ${\cal M}$
by the additional hat, $\widehat{\mathfrak{R}}$, and
by double brackets, to distinguish them from
the ordinary $\mathfrak{R}$-matrices,
acting in the representation spaces $R^{\otimes m}$.

In these terms, the statement (ii) can be given a very explicit form:
\be
\mathfrak{R}_\mu = \oplus_Q \widehat{\mathfrak{R}}_\mu^Q \otimes I_Q
\ee
i.e. the ${\mathfrak{R}}$-matrices act on the product of representations
(\ref{Rmdeco}) as
\be
{\mathfrak{R}}_\mu (R^{\otimes m}) = \oplus_Q \widehat{\mathfrak{R}}_\mu^Q ({\cal M}^Q) \otimes Q
\ee
This follows from the elementary ${\mathfrak{R}}$-matrix relation with $m=2$,
\be
{\mathfrak{R}} (R_1\otimes R_2) = \oplus_S \widehat{\mathfrak{R}}({\cal M}_{R_1R_2}^S) \otimes S
\ee
and commutativity of ${\mathfrak{R}}$-matrix and compultiplication,
\be
{\mathfrak{R}}\Delta(g){\mathfrak{R}}^{-1} = \Delta(g)
\ee

These $\widehat{\mathfrak{R}}$-matrices can be further diagonalized {\it within}
the ${\cal M}_Q$ spaces as well, but not all at once.
Each $\widehat{\mathfrak{R}}_\mu$ can be diagonalized within ${\cal M}_Q$,
but the corresponding basis
depends on $\mu$.
Basices with different $\mu$ are linearly related by
$N_{R^m}^Q\times N_{R^m}^Q$ matrices
$\widehat{\cal U}_{\,\mu\nu}^Q = (\widehat{\cal U}_{\,\nu\mu}^Q)^{-1}$.
In other words, all $\widehat{\mathfrak{R}}$-matrices can be expressed through, say, the first one,
$\widehat{\mathfrak{R}}_\mu = \widehat{\cal U}_{\mu 1}\widehat{\mathfrak{R}}_1 \widehat{\cal U}_{1\mu}$, and
\be
H_R = \sum_Q S_Q^* \ \Tr\!_{{\cal M}^Q} \left(\prod_s
\widehat{\cal U}_{\mu(s)1}\widehat{\mathfrak{R}}_1^\pm\widehat {\cal U}_{1\mu(s)} \right)
\ \ \ \ \
{\rm i.e.}
\ \ \ \ \
\K_R^Q = \Tr\!_{{\cal M}^Q}
\left(\prod_s \widehat{\cal U}_{\mu(s)1}\widehat{\mathfrak{R}}_1^\pm \widehat{\cal U}_{1\mu(s)} \right)
\label{HthrU}
\ee
where $\widehat{\mathfrak{R}}_1$ can be taken in the diagonal form.

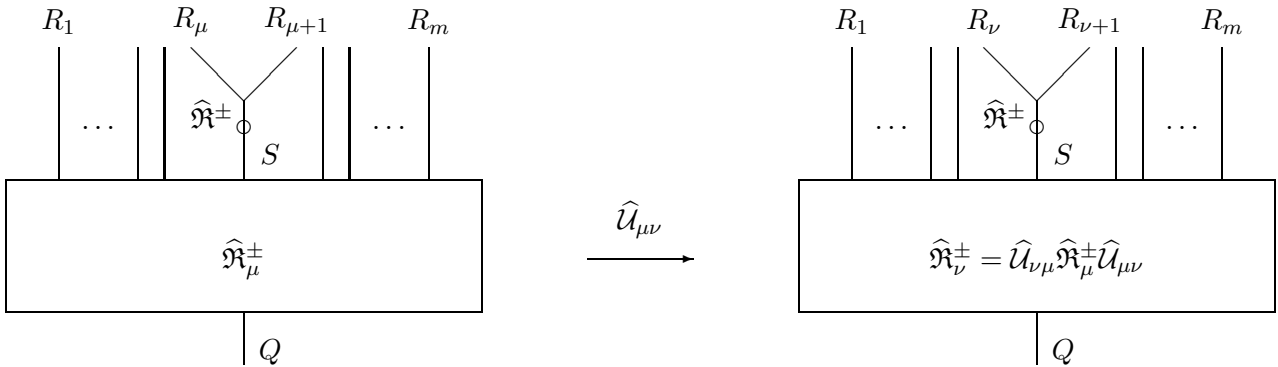
\begin{picture}(20,170)(-70,-60)
\put(0,-20){\line(0,-1){20}}
\put(0,60){\line(-1,1){20}}
\put(0,60){\line(1,1){20}}
\put(0,60){\line(0,-1){30}}
\put(-70,30){\line(0,1){50}}
\put(-40,30){\line(0,1){50}}
\put(-30,30){\line(0,1){50}}
\put(30,30){\line(0,1){50}}
\put(40,30){\line(0,1){50}}
\put(70,30){\line(0,1){50}}
\put(-90,-20){\line(0,1){50}}
\put(90,-20){\line(0,1){50}}
\put(-90,-20){\line(1,0){180}}
\put(-90,30){\line(1,0){180}}
\put(0,50){\circle{5}}
\put(-12,55){\makebox(0,0)[cc]{$\widehat{\mathfrak{R}}^{\pm}$}}
\put(-55,50){\makebox(0,0)[cc]{$\ldots$}}
\put(55,50){\makebox(0,0)[cc]{$\ldots$}}
\put(-70,90){\makebox(0,0)[cc]{$R_1$}}
\put(70,90){\makebox(0,0)[cc]{$R_m$}}
\put(-20,90){\makebox(0,0)[cc]{$R_\mu$}}
\put(20,90){\makebox(0,0)[cc]{$R_{\mu+1}$}}
\put(10,-35){\makebox(0,0)[cc]{$Q$}}
\put(10,40){\makebox(0,0)[cc]{$S$}}
\put(300,-20){\line(0,-1){20}}
\put(300,60){\line(-1,1){20}}
\put(300,60){\line(1,1){20}}
\put(300,60){\line(0,-1){30}}
\put(230,30){\line(0,1){50}}
\put(260,30){\line(0,1){50}}
\put(270,30){\line(0,1){50}}
\put(330,30){\line(0,1){50}}
\put(340,30){\line(0,1){50}}
\put(370,30){\line(0,1){50}}
\put(210,-20){\line(0,1){50}}
\put(390,-20){\line(0,1){50}}
\put(210,-20){\line(1,0){180}}
\put(210,30){\line(1,0){180}}
\put(300,50){\circle{5}}
\put(288,55){\makebox(0,0)[cc]{$\widehat{\mathfrak{R}}^{\pm}$}}
\put(245,50){\makebox(0,0)[cc]{$\ldots$}}
\put(355,50){\makebox(0,0)[cc]{$\ldots$}}
\put(230,90){\makebox(0,0)[cc]{$R_1$}}
\put(370,90){\makebox(0,0)[cc]{$R_m$}}
\put(280,90){\makebox(0,0)[cc]{$R_\nu$}}
\put(320,90){\makebox(0,0)[cc]{$R_{\nu+1}$}}
\put(310,-35){\makebox(0,0)[cc]{$Q$}}
\put(310,40){\makebox(0,0)[cc]{$S$}}
\put(130,00){\vector(1,0){40}}
\put(150,15){\makebox(0,0)[cc]{$\widehat{\cal U}_{\mu\nu}$}}
\put(0,0){\makebox(0,0)[cc]{$\widehat{\mathfrak{R}}_\mu^{\pm}$}}
\put(300,0){\makebox(0,0)[cc]{$\widehat{\mathfrak{R}}_\nu^{\pm}
= \widehat{\cal U}_{\nu\mu}\widehat{\mathfrak{R}}_\mu^{\pm}\widehat{\cal U}_{\mu\nu}$}}
\end{picture}

In the main part of the paper we are going to explain the meaning
and practical work with formula (\ref{HthrU}).
We explicitly construct the mixing matrices $\widehat{\cal U}_{1\mu}^Q$
for $R = [1]$ and $m< 5$, what allows one to express the known HOMFLY polynomials
for all the knots with up to 8 crossings through one and the same set
of $\widehat{\cal U}$-matrices.
Moreover, we construct the $\widehat{\cal U}$-matrices in a special basis,
where they can be represented as ordered products of some elementary
mixing matrices, which in the case of $[R]=1$ are essentially $2\times 2$.
The basis is labeled by different fusions of the representations $R$, i.e.
by the decorated trees inside the boxes in the picture,
and the elementary constituents of the $\widehat{\cal U}$-matrices
correspond to flipping of just one edge of the graph (therefore, their matrix elements
are given by the Racah coefficients, or fusion matrices).

The restrictions to $[R]=1$ and $m< 5$ are technical: in this case only the
Young diagrams $Q$ with no more than two columns or two rows are contributing,
and, given the $N$-independence (universality) of (\ref{HthrU}),
all $\K_R^Q$ can be derived from the representation theory of $SU_q(2)$ algebra.
When $Q$ with $l(Q)$ columns or rows appear, one needs to extend consideration
to $SU_q\left(l(Q)\right)$ at least.
This is straightforward, but in order to avoid overloading the text,
we postpone consideration of such examples to separate publications.

The explicit formulae include only group theory data,
but of two types:
some is known in general for arbitrary Young diagrams
(these are eigenvalues $r_Q$ of the quantum $R$-matrix
and values of the Schur functions $S_Q^*$, eq.(\ref{he})),
while some do not yet possess such exhaustive expressions
(these are the \LRC\ and the \CGC).

In the paper, we list the examples in the order of increasing strand number $m$.
The $2$-strand case is sensitive only to the \LRC.
The \CGC\  appear only for $3$ and more strands. We interrupt listing the examples in sect.4
with explaining the procedure general construction of the mixing matrix $U$.
After that we give the results for 4 and 5 strands. In fact, for $m=5$ and
and $R=[1]$ there is one representation, $Q=[311]$ not obtained by $SU_q(2)$
consideration.  In the present paper we restore the corresponding coefficient from the
known answers for the HOMFLY polynomials.
Some summary and comments are contained in the last section.

\section{The $2$-strand braids}

We begin with the simplest case of the two-strand braid, $m=2$:

\vspace{1cm}

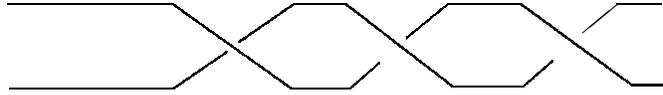
\unitlength 1mm 
\linethickness{0.4pt}
\ifx\plotpoint\undefined\newsavebox{\plotpoint}\fi 
\begin{picture}(115.5,78)(0,0)
\put(27.75,78){\line(1,0){22}}
\multiput(49.75,78)(.0471556886,-.0336826347){334}{\line(1,0){.0471556886}}
\put(65.5,66.75){\line(1,0){7.75}}
\multiput(73.25,66.75)(.038461538,.033653846){104}{\line(1,0){.038461538}}
\put(28,66.75){\line(1,0){21.75}}
\multiput(49.75,66.75)(.048657718,.033557047){149}{\line(1,0){.048657718}}
\multiput(58.5,73)(.048657718,.033557047){149}{\line(1,0){.048657718}}
\put(65.75,78){\line(1,0){7}}
\multiput(72.75,78)(.0435779817,-.0336391437){327}{\line(1,0){.0435779817}}
\put(87,67){\line(1,0){9.25}}
\multiput(96.25,67)(.038461538,.033653846){104}{\line(1,0){.038461538}}
\multiput(80.75,73.5)(.041044776,.03358209){134}{\line(1,0){.041044776}}
\put(86.25,78){\line(1,0){9.75}}
\multiput(96,78)(.0462382445,-.0336990596){319}{\line(1,0){.0462382445}}
\put(110.75,67.25){\line(1,0){4.75}}
\put(104.25,74.25){\line(6,5){4.5}}
\put(108.75,78){\line(1,0){6}}
\end{picture}

\vspace{-6cm}

\noindent
In this case we get only the torus knots and links of the type $\l[2,n]$, and
\be
H_R^{[2,n]} = \Tr_{R\otimes R}{\mathfrak{R}}^n
\ee
with integer $n$.
Knots arise for $n$ odd, links for $n$ even.
In the latter case, one can also consider
\be
H_{R_1,R_2}^{[2,n]}  = \Tr_{R_1\otimes R_2}{\mathfrak{R}}^n
\ee
with $R_1\neq R_2$.
For $R_1\otimes R_2 = \sum_Q N_{R_1R_2}^Q Q$ one gets
\be
H_{R_1,R_2}^{[2,n]}  = \sum_Q N_{R_1R_2}^Q \Tr_{Q}{\mathfrak{R}}^n
= \sum_Q N_{R_1R_2}^Q  S_Q^* r_Q^n
\label{2stbr}
\ee
where $r_Q$ is the eigenvalue of ${\mathfrak{R}}$ in the representation $Q$.
This is the final answer.

In fact, for a single ${\mathfrak{R}}$-matrix the eigenvalues $r_Q$
are known in full generality\footnote{Note that we chose in \cite{DMMSS}
the opposite sign of $\varkappa_Q$. In the case of HOMFLY polynomials, this
sing is inessential to the symmetry of ${\mathfrak{R}}$-matrix w.r.t. the
replace $q\leftrightarrow q^{-1}$. The symmetry, however, is more complicated in the
case of superpolynomials, and one has to choose the sign as in \cite{DMMSS} in that case.}:
\be
r_Q = q^{\varkappa_Q}
\label{rQ}
\ee
where $\varkappa_R$ is the eigenvalue of the simplest non-trivial
cut-and-join operator $\hat W_{[2]}$ \cite{CJ} on the character eigenfunction $S_R\{p\}$:
\be
{\rm for} \ \ \ Q = \{q_1\geq q_2\geq \ldots \geq 0\} \ \ \ \ \ \
\varkappa_Q = \frac{1}{2}\sum_i q_i(q_i+1-2i) = \nu_{Q'} - \nu_{Q},
\ \ \ \ \nu_Q =\sum_i (i-1)q_i
\ee
and $Q'$ is the transposed Young diagram.

In the simplest case of $R_1=R_2=[1]$ one has $[1]\otimes [1] = [2]+[11]$
and the two eigenvalues are just $r_{[2]} = q$ and $r_{[11]} = -\frac{1}{q}$.
Thus, (\ref{2stbr}) in this case turns into
\be
H_{[1]}^{[2,n]}
= \Tr_{[1]^{\otimes 2}} {\mathfrak{R}}^n =
\Tr_{[2]}  {\mathfrak{R}}^n + \Tr_{[11]}  {\mathfrak{R}}^n =
q^n \Tr_{[2]} I\ + \left(-\frac{1}{q}\right)^n \Tr_{[11]} I =
q^n S_2^* + \left(-\frac{1}{q}\right)^nS_{11}^*
\label{2stbrfu}
\ee

Of the knots with no more than 8 crossings, the following ones get into the set of
the $2$-strand knots:
\be
3_1 = [2,3], \ \ 5_1=[2,5], \ \ 7_1=[2,7], \ \ldots
\ee
The first notation refers, e.g., to the Rolfsen Knot Table, see \cite{katlas},
the second one is the standard (dual) notation for the torus knots.
For these knots eq.(\ref{2stbrfu}) gives:
$$
\begin{array}{c|c|cc|cc}
{\rm knot} & n & S_2^* & S_{11}^* & A & A^{-1} \\
&&&&&\\
\hline
&&&&&\\
3_1 & 3 & q^3 & -q^{-3} & q^2+q^{-2} & -1 \\
&&&&&\\
5_1 & 5 & q^5 & -q^{-5} & q^4+1+q^{-4} & -q^2-q^{-2} \\
&&&&&\\
7_1 & 7 & q^7 & -q^{-7} & q^6+q^2+q^{-2}+q^{-6} & -q^4 - 1 - q^{-4} \\
&&&&&\\
\ldots
\end{array}
$$
The two columns in the middle contain the two coefficients $h_{[1]}^{[2]}$ and $h_{[1]}^{[11]}$ which
stand in front of $S^*_{2}$ and $S^*_{11}$ correspondingly.
In this case, they are especially simple and given by the general formula (\ref{2stbrfu}).
The last two columns contain the coefficients in front of powers of $A$
in the ratio $H_{[1]}(A)/S_{[1]}^*$,
which arise from substitution of $S_{[2]}^*/S_{[1]}^*$
and $S_{[11]}^*/S_{[1]}^*$ into (\ref{2stbrfu}).
Note that, despite
\be
S_1^* = \frac{A-A^{-1}}{q-q^{-1}}, \ \ \ \
\frac{S_2^*}{S_1^*} = \frac{qA-q^{-1}A^{-1}}{q^2-q^{-2}}, \ \ \ \
\frac{S_{11}^*}{S_1^*} = \frac{q^{-1}A-qA^{-1}}{q^2-q^{-2}}
\ee
contain non-trivial denominators, they disappear from the HOMFLY polynomials.
Of course, the arising explicit expressions for the HOMFLY polynomials coincide with
the known expressions from \cite{katlas} (where $z = q-q^{-1}$, $a=1/A$, and one has
additionally to divide our expressions by the normalization factor
$A^n$).

\subsection{Colored HOMFLY polynomials}

In fact, the coefficients $\K_R^Q$ are known explicitly
in far more generality, that is,
for all torus knots $[m,n]$ \cite{knotcharexp,torus}.
In the case of $m=2$ this allows one to extend (\ref{2stbrfu}) to
arbitrary representations $R$ (i.e. to the {\it colored} HOMFLY polynomials).
For example,
\be
H_{[2]}^{[2,n]}
= \Tr_{[2]^{\otimes 2}} {\mathfrak{R}}^n = q^{6n}S^*_{4}-q^{2n}S^*_{31}+S^*_{22}    \\
H_{[11]}^{[2,n]}
= \Tr_{[11]^{\otimes 2}} {\mathfrak{R}}^n = S^*_{22}-q^{-2n}S^*_{211}+q^{-6n}S^*_{1111}    \\
H_{[3]}^{[2,n]}
= \Tr_{[3]^{\otimes 2}} {\mathfrak{R}}^n = q^{15n}S^*_{6}-q^{9n}S^*_{51}+q^{5n}S^*_{42}-q^{3n}S^*_{33}    \\
H_{[21]}^{[2,n]}
= \Tr_{[21]^{\otimes 2}} {\mathfrak{R}}^n = q^{5n}S^*_{42}-q^{3n}S^*_{411}-q^{3n}S^*_{33}+
q^{-3n}S^*_{3111}+q^{-3n}S^*_{222}-q^{-5n}S^*_{2211}    \\
H_{[111]}^{[2,n]}
= \Tr_{[111]^{\otimes 2}} {\mathfrak{R}}^n = q^{-3n}S^*_{222}-q^{-5n}S^*_{2211}+q^{-9n}S^*_{21111}
-q^{-15n}S^*_{111111}    \\
\ldots
\ee

\subsection{$2$-strand superpolynomials}

According to \cite{DMMSS}, the switch from HOMFLY to superpolynomials
in the case of toric knots is straightforward.
In the $2$-strand case
it is enough to substitute the two ${\mathfrak{R}}$-matrix eigenvalues
\be
q,\ -\frac{1}{q} \ \longrightarrow \ q,\ -\frac{1}{t}
\ee
and the two quantum dimensions
\be
\frac{S_2^*}{S_1^*} \ \longrightarrow \frac{M_2^*}{M_1^*} = {\{Aq\}\over\{qt\}}\ ,
\hspace{2cm}
\frac{S_{11}^*}{S_1^*} \ \longrightarrow \frac{1-t^4}{1-q^2t^2}\frac{M_{11}^*}{M_1^*} = {\{At^{-1}\}\over\{t^2\}}
\ee
The result is
\be
P_{[1]}^{[2,2k+1]} = q^{2k+1}M_2^* - \frac{1-t^4}{1-q^2t^2}\ {q\over t}\ t^{-2k-1}M_{11}^*
= -\frac{M_1^*}{A}{q\over t^{2k+1}} \left[ {1-(q^2t^2)^{k+1}\over 1-q^2t^2}(-A^2)     +
t^2{1-(q^2t^2)^k\over 1-q^2t^2} \right]
\ee
All the coefficients in the emerging polynomial (in the square brackets)
are positive integers
which are related to dimensions of the Khovanov-Rozhansky homologies \cite{KhR,DGR}.

\section{The $3$-strand braids}

This is the first case where there are at least two different ${\mathfrak{R}}$-matrices,
${\mathfrak{R}}_1$ and ${\mathfrak{R}}_2$, and a mixing matrix $U_{12}$ arise for the first time.
On the other hand, since the decomposition
\be
[1]^3 = [3] + 2\,[21] + [111]
\label{1cube}
\ee
contains exactly the same number of different representations (three)
as there are different powers of $A^2$ in $S^*$ at the third level,
the character expansion of the HOMFLY polynomials in this case is defined unambiguously
and can be checked in an independent way.

Since there are $N_{1^2}^{[21]}=2$ representations $[21]$ in the decomposition (\ref{1cube}),
the $U$-matrix will be non-trivial only in this sector, $Q=[21]$,
and it will be a $2\times 2$ orthogonal matrix,
\be
U = \ldb\begin{array}{cc} C & S \\ -S & C\end{array}\rdb\ ,
\hspace{1cm}
\hbox{while}
\hspace{2cm}
{\mathfrak{R}}_1 = {\mathfrak{R}}\otimes I
\label{2sU}
\ee
in this sector is just
${\mathfrak{R}}_1 = \ldb\begin{array}{cc} q & 0 \\ 0 & -\frac{1}{q}\end{array}\rdb$.
We remind that the ${\mathfrak{R}}$-matrices in the ${\cal M}$ space are denoted by the double brackets,
in order to distinguish them from the ordinary ${\mathfrak{R}}$-matrices,
acting in the space of $[1]^3$.
The ordinary ${\mathfrak{R}}$ matrices depend on $N$ (are of the size $N^3\times N^3$,
while the double bracket (or hatted) ones are $1\times 1$, $2\times 2$ and $1\times 1$
in the sectors $Q=[3]$, $Q=[21]$ and $Q=[111]$ respectively.

An arbitrary $3$-strand braid is parameterized by a sequence of integers
$a_1,b_1,a_2,b_2,\ldots$ (in this figure $a_1=-2$, $b_1=2$, $a_2=-1$, $b_2=3$:
this is knot $8_{10}$ ):

\vspace{0.7cm}

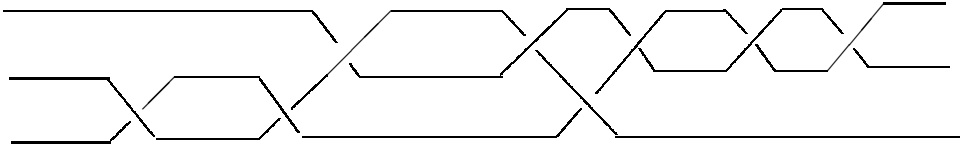
\unitlength 1mm 
\linethickness{0.4pt}
\ifx\plotpoint\undefined\newsavebox{\plotpoint}\fi 
\begin{picture}(145.5,53)(0,0)
\put(19.5,34.5){\line(1,0){13.25}}
\put(41.25,43.25){\line(1,0){11.25}}
\put(19.25,43){\line(1,0){13.25}}
\put(38.75,35){\line(1,0){13.75}}
\put(61.25,43.25){\line(1,1){8.75}}
\put(70,52){\line(1,0){14.75}}
\put(18.5,52){\line(1,0){41}}
\multiput(59.5,52)(.033505155,-.043814433){97}{\line(0,-1){.043814433}}
\put(58.25,35.25){\line(1,0){33.75}}
\multiput(92,35.25)(.033505155,.038659794){97}{\line(0,1){.038659794}}
\multiput(64.5,45)(.03289474,-.04605263){38}{\line(0,-1){.04605263}}
\put(65.75,43.25){\line(1,0){19}}
\multiput(84.5,43.5)(.0346153846,.0336538462){260}{\line(1,0){.0346153846}}
\multiput(84.75,52)(.03370787,-.03651685){89}{\line(0,-1){.03651685}}
\multiput(52.5,43)(.033653846,-.046474359){156}{\line(0,-1){.046474359}}
\multiput(52.5,35)(.03353659,.03353659){82}{\line(0,1){.03353659}}
\multiput(56.75,39)(.035447761,.03358209){134}{\line(1,0){.035447761}}
\multiput(32.25,43)(.033602151,-.041666667){186}{\line(0,-1){.041666667}}
\multiput(32.75,34.75)(.03333333,.03333333){75}{\line(0,1){.03333333}}
\put(37,39){\line(1,1){4.25}}
\put(99.75,35.25){\line(1,0){45.75}}
\multiput(100,35.5)(-.0336990596,.0352664577){319}{\line(0,1){.0352664577}}
\multiput(97.25,41)(.0336363636,.04){275}{\line(0,1){.04}}
\put(106.5,52){\line(1,0){7.75}}
\put(121.25,44){\line(1,0){6.75}}
\put(128,44){\line(5,6){7.5}}
\put(135.5,53){\line(1,0){8.25}}
\put(93.25,52.25){\line(1,0){5.75}}
\multiput(99,52.25)(.03353659,-.04268293){82}{\line(0,-1){.04268293}}
\multiput(103,47)(.03333333,-.05){60}{\line(0,-1){.05}}
\put(105,44){\line(0,1){0}}
\put(105,44){\line(1,0){9.5}}
\multiput(114.5,44)(.033632287,.036995516){223}{\line(0,1){.036995516}}
\put(122,52.25){\line(1,0){5.25}}
\multiput(127.25,52.25)(.03353659,-.03963415){82}{\line(0,-1){.03963415}}
\multiput(131.5,47)(.03333333,-.04166667){60}{\line(0,-1){.04166667}}
\put(133.5,44.5){\line(1,0){10.75}}
\multiput(114.25,52.25)(.03370787,-.03651685){89}{\line(0,-1){.03651685}}
\multiput(121,44)(-.03333333,.04666667){75}{\line(0,1){.04666667}}
\end{picture}

\vspace{-2.7cm}

\noindent
The corresponding $H_{[1]}$ is given by
\be
H_{[1]}=
\Tr \left\{\Big({\mathfrak{R}}\otimes I\Big)^{a_1} \Big(I \otimes {\mathfrak{R}}\Big)^{b_1}
\Big({\mathfrak{R}}\otimes I\Big)^{a_2} \Big(I \otimes {\mathfrak{R}}\Big)^{b_2} \ldots\right\} \
\\
=\sum_{Q=[111],[21],[3]}\tr\left\{\left(\hat{\mathfrak{R}}_1^Q\right)^{a_1}
\left(\hat{\mathfrak{R}}_2^Q\right)^{b_1}\left(\hat{\mathfrak{R}}_1^Q\right)^{a_2}
\left(\hat{\mathfrak{R}}_2^Q\right)^{b_2}\ldots\right\}=\nn
\ee
\fr{
=q^{a_1+b_1+a_2+b_2+\ldots} S_3^* +
\left(-\frac{1}{q}\right)^{a_1+b_1+a_2+b_2+\ldots}S_{111}^* +\hspace{2cm}
  \\
+ \tr \left\{
\ldb\begin{array}{cc} q & 0 \\ 0 & -\frac{1}{q}\end{array}\rdb^{a_1}
\underbrace{\ldb\begin{array}{cc} C & S \\ -S & C\end{array}\rdb
\ldb\begin{array}{cc} q & 0 \\ 0 & -\frac{1}{q}\end{array}\rdb^{b_1}
\ldb\begin{array}{cc} C & -S \\ S & C\end{array}\rdb}
\times
\right. \\  \left. \times
\ldb\begin{array}{cc} q & 0 \\ 0 & -\frac{1}{q}\end{array}\rdb^{a_2}
\underbrace{\ldb\begin{array}{cc} C & S \\ -S & C\end{array}\rdb
\ldb\begin{array}{cc} q & 0 \\ 0 & -\frac{1}{q}\end{array}\rdb^{b_2}
\ldb\begin{array}{cc} C & -S \\ S & C\end{array}\rdb}
\ldots \right\}
S_{21}^*}
We took the orthogonality of $U$ into account.
We evaluate $C$ and $S$ explicitly in s.4 below,
it turns out that
\be
C = \frac{1}{[2]_q} = \frac{1}{q+q^{-1}}, \ \ \ \ \
S = \frac{\sqrt{[3]_q}}{[2]_q} =  \frac{\sqrt{q^2+1+q^{-2}}}{q+q^{-1}}
\ee
In the classical limit of $q\to 1$ this orthogonal matrix $U$ describes rotating
at angle $\pi/3$, while at generic $q$ the rotation angle $\theta$ ($C=\cos\theta$,
$S=\sin\theta$) is not that nice.
{\bf Note that {\it all} the $3$-strand knots are described by a single formula
with just two non-trivial entries (\ref{2sU}).}
All the corresponding HOMFLY polynomials from \cite{katlas} are certainly immediately
reproduced
(we remind that in \cite{katlas} $z = q-q^{-1}$ and $a=1/A$), see the Table at the next page.

\paragraph{Comment on the Table.}
$3_1$ is a torus knot $[2,3]$, therefore, it has a natural minimal braid
realization as a $2$-strand braid, and in this quality it was already considered
in the previous section. However, since $[2,3]=[3,2]$, it also possesses a $3$-strand representation
and, hence, is also present in the Table.
The corresponding ${\cal H}_R\{p\}\equiv\sum_R h_R^QS_Q(p)$ are, of course, different:
\be
{\cal H}_R^{[3,2]}(p) \neq {\cal H}_R^{[2,3]}(p)
\ee
i.e. ${\cal H}_R(p)$'s are {\it not} knot invariants,
but their restrictions on the subspace {\it do} coincide:
\be
H_R^{[3,2]} = {\cal H}_R^{[3,2]}\{p^*\} = {\cal H}_R^{[2,3]}\{p^*\} = H_R^{[2,3]}
\ee
and this is indeed seen in the last columns of the Tables of this and the previous sections (for $R=[1]$).

$8_{19}$ is a torus knot $[3,4]$, therefore, it possesses also another braid representation,
but this time it is also a $3$-strand braid.
Because of this, for these two representation not only $H_R^{8_{19}}$ are the same,
but also the entire ${\cal H}_R^{8_{19}}$, as follows from the Table (for $R=[1]$):
the coefficients $\K_R^Q$ are the same for both realizations.

We added two more lines to the Table, describing knot $10_{139}$.
The task is to explain the coincidence observed in \cite{DMMSS}.
$10_{139}$ can be considered as $5_2$ with an attached torus braid:
\be
{\cal B}^{10_{139}} = {\cal B}^{5_2} \Big(({\mathfrak{R}}\otimes I)(I\otimes {\mathfrak{R}})\Big)^3
\label{10139}
\ee
and, thus, its HOMFLY polynomial (and even the superpolynomial) can be obtained from
that of $5_2$ by a "torus evolution" described in \cite{DMMSS}.
In our Table we demonstrate that the HOMFLY polynomial for (\ref{10139}) is indeed the same
as for the conventional realization of $10_{139}$.

\newpage

\thispagestyle{empty}

\hspace{-1cm}
{\tiny{\rotate{
\begin{tabular}{|c|c|ccc|ccc|}
\multicolumn{8}{c}{{\small\bf The Table of HOMFLY polynomials including 3-strand
knots with no more than 8 crossings plus knot $10_{139}$ (see (\ref{10139}))}}\\
\multicolumn{8}{c}{}\\
\hline
&&&&&&&\\
{\rm knot}& $(a_1,b_1,\ldots)$&
$S_3^*$ & $S_{21}^*$ & $S_{111}^*$ & $A^2$ & $A^0$ & $A^{-2}$ \\
&&&&&&&\\
\hline
&&&&&&&\\
$3_1$& (-1,-1,-1,-1)&$q^{-4}$ &-1&$q^4$&0&-1&$q^{2}+q^{-2}$\\
&&&&&&&\\
$4_1$& (1,-1,1,-1)&1&$q^4-2q^2+1-2q^{-2}+q^{-4}$&1&1&$-q^{2}+1-q^{-2}$&1\\
&&&&&&&\\
$5_2$& (-1,1,-1,-3)&$q^{-4}$&$-(q^4-q^2+1-q^{-2}+q^{-4})$&$q^{4}$&-1&$q^{2}-1+q^{-2}$&$q^{2}-1+q^{-2}$\\
&&&&&&&\\
$6_2$& (1,-1,1,-3)&$q^{-2}$&$q^6-2q^4+2q^2-3+2q^{-2}-2q^{-4}+q^{-6}$&$q^2$&$q^{2}-1+q^{-2}$&$-q^{4}+q^{2}-2+q^{-2}-q^{-4}$
&$q^{2}+q^{-2}$\\
&&&&&&&\\
$6_3$& (2,-1,1,-2)&1&$-(q^2-1+q^{-2})(q^4-q^2+1-q^{-2}+q^{-4})$&1&$-q^{2}+1-q^{-2}$&$q^{4}-q^{2}+3-q^{-2}+q^{-4}$&
$-q^{2}+1-q^{-2}$\\
&&&&&&&\\
$7_3$& (1,-1,1,5)&$q^6$&$-q^6+q^4-2q^2+3-2q^{-2}+q^{-4}-q^{-6}$&$q^{-6}$&$q^{4}-q^{2}+1-q^{-2}+q^{-4}$&
$q^{4}-q^{2}+2+q^{-4}-q^{-2}$&$-q^{2}-q^{-2}$\\
&&&&&&&\\
$7_5$& (-2,1,-1,-4)&$q^6$&$-(q^2-1+q^{-2})(q^4-q^2+1-q^2+1)$&$q^6$
&$-q^{2}+1-q^{-2}$&$q^{4}-2q^{2}+2-2q^{-2}+q^{-4}$&$q^{4}-q^{2}+2+q^{-4}-q^{-2}$\\
&&&&&&&\\
$8_2$& (1,-1,1,-5)&$q^{-4}$&$q^{8}-2q^{6}+2q^{4}-3q^{2}+3-3q^{-2}+2q^{-4}-2q^{-6}+q^{-8}$&$q^4$
&$q^{4}-q^{2}+1-q^{-2}+q^{-4}$&$-q^{6}+q^{4}-2q^{2}+1-2q^{-2}+q^{-4}-q^{-6}$&$q^{4}+1+q^{-4}$\\
&&&&&&&\\
$8_5$& (-1,3,-1,3)&$q^4$&$q^{8}-2q^{6}+3q^{4}-4q^{2}+3-4q^{-2}+3q^{-4}-2q^{-6}+q^{-8}$&$q^{-4}$
&$q^{4}+2+q^{-4}$&$-q^{6}+q^{4}-3q^{2}+1-3q^{-2}+q^{-4}-q^{-6}$&$q^{4}-q^{2}+2+q^{-4}-q^{-2}$\\
&&&&&&&\\
$8_7$& (-2,1,-1,4)&$q^2$&$-q^{8}+2q^{6}-3q^{4}+4q^{2}-5+4q^{-2}-3q^{-4}+2q^{-6}-q^{-8}$&$q^{-2}$
&$-q^{4}+q^{2}-1+q^{-2}-q^{-4}$&$q^{6}-q^{4}+3q^{2}-2+3q^{-2}-q^{-4}+q^{-6}$&$-q^{4}+q^{2}-2+q^{-2}-q^{-4}$\\
&&&&&&&\\
$8_9$& (3,-1,1,-3)&1&$(q^2-1+q^{-2})(q^6-q^4+q^2-3+q^{-2}-q^{-4}+q^{-6})$&1
&$q^{4}-q^{2}+2+q^{-4}-q^{-2}$&$-q^{6}+q^{4}-3q^{2}+3-3q^{-2}+q^{-4}-q^{-6}$&$q^{4}-q^{2}+2+q^{-4}-q^{-2}$\\
&&&&&&&\\
$8_{10}$& (-2,2,-1,3)&$q^2$&$-q^{8}+2q^{6}-4q^{4}+5q^{2}-5+5q^{-2}-4q^{-4}+2q^{-6}-q^{-8}$&$q^{-2}$
&$-q^{4}+q^{2}-2+q^{-2}-q^{-4}$&$q^{6}-q^{4}+4q^{2}-2+4q^{-2}-q^{-4}+q^{-6}$&$-q^{4}+q^{2}-3+q^{-2}-q^{-4}$\\
&&&&&&&\\
$8_{16}$& (1,-1,1,-2,1,-2)&$q^{-2}$&$-q^{8}+3q^{6}-5q^{4}+6q^{2}-7+6q^{-2}-5q^{-4}+3q^{-6}-q^{-8}$&$q^2$
&$-q^{4}+2q^{2}-3+2q^{-2}-q^{-4}$&$q^{6}-2q^{4}+4q^{2}-4+4q^{-2}-2q^{-4}+q^{-6}$&$-q^{4}+2q^{2}-2+2q^{-2}-q^{-4}$\\
&&&&&&&\\
$8_{17}$& (2,-1,1,-1,1,-2)&1&$q^{8}-3q^{6}+5q^{4}-7q^{2}+7-7q^{-2}+5q^{-4}-3q^{-6}+q^{-8}$&1
&$q^{4}-2q^{2}+3-2q^{-2}+q^{-4}$&$-q^{6}+2q^{4}-4q^{2}+5-4q^{-2}+2q^{-4}-q^{-6}$&$q^{4}-2q^{2}+3-2q^{-2}+q^{-4}$\\
&&&&&&&\\
$8_{18}$& (1,-1,1,-1,1,-1,1,-1)&1&$q^{8}-4q^{6}+6q^{4}-8q^{2}+9-8q^{-2}+6q^{-4}-4q^{-6}+q^{-8}$&1
&$q^{4}-3q^{2}+3-3q^{-2}+q^{-4}$&$-q^{6}+3q^{4}-4q^{2}+7-4q^{-2}+3q^{-4}-q^{-6}$&$q^{4}-3q^{2}+3-3q^{-2}+q^{-4}$\\
&&&&&&&\\
$8_{19}$& (1,3,1,3)=[3,4]&$q^8$&-1&$q^{-8}$
&$q^{6}+q^{2}+1+q^{-2}+q^{-6}$&$-q^{4}-q^{2}-1-q^{-2}-q^{-4}$&1\\
& (1,1,1,1,1,1,1,1)&$q^8$ &-1& $q^{-8}$
&$q^{6}+q^{2}+1+q^{-2}+q^{-6}$&$-q^{4}-q^{2}-1-q^{-2}-q^{-4}$&1\\
&&&&&&&\\
$8_{20}$& (-1,-3,-1,3)&$q^{-2}$&$-q^{6}+q^{4}-q^{2}+1-q^{-2}+q^{-4}-q^{-6}$&$q^2$
&$-q^{2}-q^{-2}$&$q^{4}+2+q^{-4}$&$-q^{2}+1-q^{-2}$\\
&&&&&&&\\
$8_{21}$& (-2,2,-1,-3)&$q^{-4}$&$q^{6}-2q^{4}+2q^{2}-3+2q^{-2}-2q^{-4}+q^{-6}$&$q^4$
&$q^{2}-1+q^{-2}$&$-q^{4}+q^{2}-3+q^{-2}-q^{-4}$&$2q^{2}-1+2q^{-2}$\\
&&&&&&&\\
&&&&&&&\\
\hline
&&&&&&&\\
$10_{139}$& (2,3,1,4)&$q^{10}$&$-\frac{q^5+q^{-5}}{q+q^{-1}}$&$q^{-10}$
&$q^{8}+q^{4}+q^{2}+q^{-2}+q^{-4}+q^{-8}$&$-q^{6}-q^{4}-2-q^{-4}-q^{-6}$&$q^{2}-1+q^{-2}$\\
& (1,-1,1,3,1,1,1,1,1,1)&$q^{10}$&$-\frac{q^5+q^{-5}}{q+q^{-1}}$&$q^{-10}$
&$q^{8}+q^{4}+q^{2}+q^{-2}+q^{-4}+q^{-8}$&$-q^{6}-q^{4}-2-q^{-4}-q^{-6}$&$q^{2}-1+q^{-2}$\\
&&&&&&&\\
&&&&&&&\\
&&&&&&&\\
\hline
\end{tabular}
}}}

\subsection*{Comments}

\subsubsection*{1. Torus knots}

Among arbitrary $3$-strand braids described by arbitrary sequences $a_1,b_1,\ldots$,
one can select particular {\it series}.

The simplest example is the series of {\it torus} knots $[3,n]$ (for $n=3k$ these are links),
with $a_k=b_k=1$, $k=1,\ldots,n$.
In this case it makes sense to diagonalize not ${\mathfrak{R}}_1$,
but the product ${\mathfrak{R}}_1{\mathfrak{R}}_2$, and the corresponding eigenvalues
$r_Q^{torus}(i)$, $i = 1,\ldots,{\rm dim} {\cal M}_Q = N_{1^2}^Q$
define the generic expression
\be
H_{[1]}^{[3,n]} = \sum_Q \sum_{i=1}^{N_{1^2}^Q} \left(r_Q^{torus}(i)\right)^n S_Q^*
\ee
Since
\be
\Big({\mathfrak{R}}\otimes I\Big)\Big( I\otimes {\mathfrak{R}}\Big)
= \ldb \begin{array}{cccc}
q^2&&& \\ & -\frac{1}{q[2]_q} & \frac{q\sqrt{[3]_q}}{[2]_q} & \\
& -\frac{\sqrt{[3]_q}}{q[2]_q} & -\frac{q}{[2]_q} & \\
&&& \left(-\frac{1}{q}\right)^2
\end{array}\rdb,
\ee
where the central $2\times 2$-block is just
\be
\left(\begin{array}{cc}
q&\\
&-{1\over q}
\end{array}\right)\cdot\left(\begin{array}{cc}
c&s\\
-s&c
\end{array}\right)\cdot\left(\begin{array}{cc}
q&\\
&-{1\over q}
\end{array}\right)\cdot\left(\begin{array}{cc}
c&-s\\
s&c
\end{array}\right)=\left(\begin{array}{cc}
-\frac{1}{q[2]_q} & \frac{q\sqrt{[3]_q}}{[2]_q} \\
-\frac{\sqrt{[3]_q}}{q[2]_q} & -\frac{q}{[2]_q}
\end{array}\right)
\ee
one obtains
\be
r_{[3]}^{torus} = q^2,  \\
r_{[21]}^{torus}(1) = e^{2\pi i/3},\ \ \  \ \ \ \ \ \ \
r_{[21]}^{torus}(2) = e^{-2\pi i/3},  \\
r_{[111]}^{torus} = -\frac{1}{q^2},  \\
\ee
i.e.
\be
\sum_{i=1}^{N_{1^2}^{[21]}} \left(r_{[21]}^{torus}(i)\right)^n
 =\left(e^{2\pi i\over 3}\right)^n+\left(e^{-2\pi i\over 3}\right)^n
 =2\cos{2\pi n\over 3}= \left\{ \begin{array}{ccc}
 -1 & {\rm for} & n=2k\pm 1, \\
 2 & {\rm for} & n =2k
 \end{array}\right.
\ee
Thus, finally,
\be
H_R^{[3,n]} = q^{2n}S^*_3+2\cos {2\pi n\over 3}\cdot S^*_{21}+q^{-2n}S^*_{111}
\ee

\subsubsection*{2. Other series}

Of course, one can consider not only the series of torus knots, but any other series.

For example, take $a_k=1$, $b_k=3$, $k=1,\ldots,n$.
In this case in the $[21]$ sector one has
\be
\Big({\mathfrak{R}}\otimes I\Big)\Big( I\otimes {\mathfrak{R}}\Big)^3
= \ldb \begin{array}{cccc}
q^{4}&&& \\ &\frac{q^6-q^4-1}{q^3[2]}&\frac{(q^4-q^2-1)\sqrt{[3]_q}}{q[2]_q}&\\
&-\frac{(q^4-q^2-1)\sqrt{[3]_q}}{q^3[2]_q}&-\frac{q^6-q^4-1}{q^3[2]_q}& \\
&&& \frac{1}{q^{4}}
\end{array}\rdb
\ee
with eigenvalues
\be
e^{\pm 2\pi i/3}
\ee
just the same as for the torus case.
However, in the sectors $[3]$ and $[111]$ we now have eigenvalues
$q^{4n}$ and $q^{-4n}$, instead of $q^{2n}$ and $q^{-2n}$ in the torus case,
thus, expressions for the HOMFLY polynomials will be different.

Likewise, for the similar series $a_k=a$, $b_k=b$, $k=1,\ldots,n$ with other values of $a$ and $b$
one has $q^{(a+b)n}$ and $(-1/q)^{(a+b)n}$ in the $[3]$ and $[111]$ sectors respectively,
while in the $[21]$ sector the situation is more complicated,
the two eigenvalues being

\be
r_{[21]}^{(a,b)}(\pm) = \frac{1}{2}\alpha\pm\frac{1}{2}\sqrt{\alpha^2-4(-1)^{a+b}}
\\ \nonumber
\alpha=\frac {q^{a+b}+(-1)^{a+b}q^{-a-b}+[3]_q\left((-1)^bq^{a-b}+
(-1)^aq^{b-a}\right)}{[2]_q^2}
\ee
and only the first few eigenvalues are simple and do not contain square roots:
$r_{[21]}^{(1,1)}=r_{[21]}^{(1,3)}=r_{[21]}^{(3,1)}=
\displaystyle{e^{\pm {2\pi i\over 3}}}$, $r_{[21]}^{(1,2)}=r_{[21]}^{(2,1)}=\pm 1$
and $r_{[21]}^{(2,2)}=q^{\pm 2}$.

Of course, the square roots disappear from the sums $\Big(r_{[21]}(+)\Big)^n +
\Big(r_{[21]}(-)\Big)^n$.

In deserves noting that the triviality of knots with $(a,b)$ equal to $(0,0)$ and
$(1,0)$ implies respectively the following identities between the Schur functions
\be
\Tr_{[21]} I = 2 \ \Longrightarrow\  (S_3+2S_{21}+S_1^3) = S_1^3 ,   \\
\Tr_{[21]} {\mathfrak{R}} = q-\frac{1}{q} \ \Longrightarrow \ \
\left(qS_3^* + \left(q-\frac{1}{q}\right)S_{21}^* -\frac{1}{q}S_{111}^*\right) =
S_1^*\left(q S_2^* - \frac{1}{q}S_{11}^*\right),  \\
\ldots
\ee

\subsubsection*{3. Composite knots}

If one takes arbitrary numbers of crossings, one can get not only prime knots
which can not be represented as a sum of two or more independent knots
(successively made one after another on the same string) or
links
from the tables. There is also another possibility: it can be composite of knots or
links. For each combination of numbers of crossings one can evaluate the HOMFLY
polynomial using our method (see Table below for examples). It is quite easy
to distinguish these cases by the form of the HOMFLY polynomial. If it has several
multipliers then it is a composite link or knot, and each of the multipliers is a
HOMFLY polynomial of the corresponding knot or link (see (\ref{redkn}) for examples).
If it is a link then for each additional loop it has a multiplier
$\frac{1}{(q-q^{-1})}$. Also links and knots can be oriented differently. For knots
the change of the orientation is a quite simple procedure: one just should take the
opposite to all the numbers of the crossings and it corresponds to the
substitution of $A^{-1}$ instead of $A$ into the HOMFLY polynomial.
With the links, however, the situation is more complicated, because one can change the
direction of some of the loops. In this case, the HOMFLY polynomial changes non-trivially
(see (\ref{orunkn}) for example: these two links have the same picture, but the
different orientation), but the change of numbers of crossings is still quite trivial:
one should change all the numbers which the changed loop involves.

The HOMFLY polynomials of the composite knots are always factorized:
\be
H_R^{{\cal K}_1\sharp {\cal K}_2} = \frac{H_R^{{\cal K}_1}
H_R^{{\cal K}_2}}{H_R^\emptyset}
\label{Hredeco}
\ee
This is slightly different from the decomposition of the HOMFLY polynomials for
the disjoint unification of two knots,
\be
H_R^{{\cal K}_1\cup {\cal K}_2} = H_R^{{\cal K}_1} H_R^{{\cal K}_2}
\ee

The simplest examples of composite knots are made out of the $3$-strand
braids.
In particular, composite are obviously the knots with two non-vanishing odd parameters
$a_1$ and $b_1$, they are in fact a composition
of two $2$-strand knots with parameters $a_1$ and $b_1$ respectively. Similarly,
$\{a_1,b_1,c_1\}$ and $\{a_1,b_1,c_1,d_1\}$ will be obviously composite knots
(or links) for 4 and 5 strands respectively.

Of course, the set of composite $3$-strand knots is not exhausted by
the set $(a_1,b_1)$.
In particular, the "torus descendant" of $4_1$,
which was considered in \cite{DMMSS},
in our present notation it is a $3$-strand braid $(1,-1,1,-1,1,1,1,1,1,1)$,
is equivalent to $(3,3)$:
\be
H_{[1]}^{(1,-1,1,-1,1,1,1,1,1,1)} = H_{[1]}^{(3,3)} =
\Big(A(q^2+q^{-2}) - A^{-1}\Big)^2 S_1^*  =
\frac{ \left(H_{[1]}^{3_1}\right)^2}{S_1^*}
\ee

These are examples of the simplest composite knots (the sequences $a_1,b_1,\ldots$
are written here in braces, $\#$ denotes composition of knots and we use dual notations
for the torus knots):

\be
\label{redkn}
\text{Knot}=\{3,3\}=T[2,3]\#T[2,3]:\hspace{4cm}
H=\left(A(q^2+q^{-2})-A^{-1}\right)^2
\ee
\be
\text{Knot}=\{3,4\}=T[2,3]\#T[2,4]:\hspace{10cm}\\
\hspace{2cm}
H=\frac{1}{q-q^{-1}}\left(A(q^2+q^{-2})-A^{-1}\right)\left(A(q^4-q^2+1-
q^{-2}+q^{-4})-A^{-1}(q^2-1+q^{-2})\right)
\ee

\be
\label{orunkn}
\text{Knot}=\{2,1,1,1\}=T[2,4]:\hspace{1.5cm}
H=\frac{1}{q-q^{-1}}\left(A^2(q^4-q^2+1-q^{-2}+q^{-4})-(q^2-1+q^{-2})\right) \\
\text{Knot}=\{2,1,-1,1\}=L4a1:\hspace{2cm}
H=\frac{1}{q-q^{-1}}\left(A^2(q^2-2+q^{-2})+(q^2-1+q^{-2})-A^{-2}\right)
\ee
and
\be
H^{T[2,3]}=A(q^2+q^{-2})-A^{-1}\\
H^{T[2,4]}=\frac{1}{q-q^{-1}}\left(A(q^4-q^2+1-q^{-2}+q^{-4})-A^{-1}(q^2-1+q^{-2})\right)
\ee

More examples of knots/links
with $a_1,b_1,a_2,b_2$ not exceeding 2
in absolute value can be found in the following Table:

\bigskip

\label{12crkn}
\begin{tabular}{c|c|c|c}
knots & & & \\
\hline
$0_1$& 1,-1,-1,-1 & 1,1,-1,-1 & 1,1,-2,2\\
$3_1$&-1,-1,-1,-1&1,1,-2,-2& \\
$4_1$&1,-1,1,-1&1,-1,-2,2& \\
$5_1$&1,1,2,2&&\\
$5_2$&1,-1,-2,-2&&\\
$6_3$&1,-2,2,-1&&\\
\hline
links & & & \\
$2\times 0_1$&1,1,1,-2&&\\
$L2a1$&1,-1,-1,-2&1,1,-1,-2&\\
$L4a1$&1,-1,-2,-1&2,-2,-1,-2&\\
$T(2,4)$&2,1,1,1&&\\
$L5a1$&1,-1,1,-2&1,2,-2,-2&\\
$L6n1$&1,2,1,-2&&\\
$L7n1$&-1,-2,-2,-2&&\\
$L7n2$&2,-1,-2,-2&&\\
$4L2a1\#L2a1$&1,2,-1,2&&\\
$L2a1\#L2a1'$&1,2,-1,-2&&\\
\hline
\end{tabular}

\bigskip

A large variety of links and composite knots
explains why the sequences $a_1,b_1,\ldots$ appearing in our Table
of 3-strand knots above
are not arbitrary.
The list actually refers to the prime knots.

As follows from the observation in \cite{DMMSS}
on the "descendant"
$(1,-1|1,-1|1,1|1,1|1,1)$ of the knot $4_3 = (1,-1|1,-1)$,
which (descendant) appears to be a composite knot,
the "{\it superpolynomials}" for composite knots are more
complicated than those for the prime ones.
This is not very surprising because of the special role of
the unknot $H_R^\emptyset$  in the decomposition rule (\ref{Hredeco}),
and unknot superpolynomial is a subtle issue
already by itself.

\section{Systematic description of ${\cal U}$-matrices}

Before proceeding to the $4$-strand and $5$-strand knots,
we now provide a general construction of the mixing $U$-matrices.
In fact, as we already noted they are nothing but the Racah coefficients, however,
as we shall see, there might be recovered a deep hierarchical structure expressing the
relevant compositions of the Racah coefficients through simpler ingredients,
which seem to possess a general (universal) description.

In application to braid calculus it will be convenient to
build ${\cal U}$-matrices hierarchically, and the fusion pattern of representations
can be presented by rooted
tree diagrams.
The key point is that the elementary ${\mathfrak{R}}$-matrices
appear only at the first level of hierarchy,
when just {\it two} representations $R$ merge, while
mixing arises at different levels of the tree. We count levels
from the top, not from the root.
All our mixings are actually described by $2\times 2$ matrices.
Moreover, they are universal in the following sense:

\begin{itemize}
\item
everything unmixed at one level remains unmixed at the next levels
\item
all descendants of a given mixing are the same at all next levels
\end{itemize}

Then the full ${\cal U}$-matrices are represented as the ordered products
of elementary ones, arising at different levels.
Our simple examples (for $m\leq 5$) involve just three levels,
and we denote the corresponding elementary factors by $U$, $V$ and $W$.
For $R=[1]$ all the three will be essentially $2\times 2$ matrices.
Of course, one can ignore this additional hierarchical structure
and evaluate the full ${\cal U}$ matrices directly.

It is sufficient to deal with $SU_q(2)$ representations in order to study
all the representations given by the Young diagrams with no more than two rows.
After that one can use the symmetry $R,q \longrightarrow R',-1/q$
in order to find expressions for the two-column representations.
For $[1]^{\otimes m}$ the $SU_q(2)$ is insufficient starting from $m=5$,
when $[311]$ arises.

We study the mixing, looking explicitly at the highest weights of
the corresponding representations.
The highest weights, in turn, show up at different {\it grades} of
the Verma module constructions.

Let us look at the example of tensor product of four fundamental representations
$[1]^{\otimes 4}$. It can be expanded, for instance, in the following way:

\bigskip

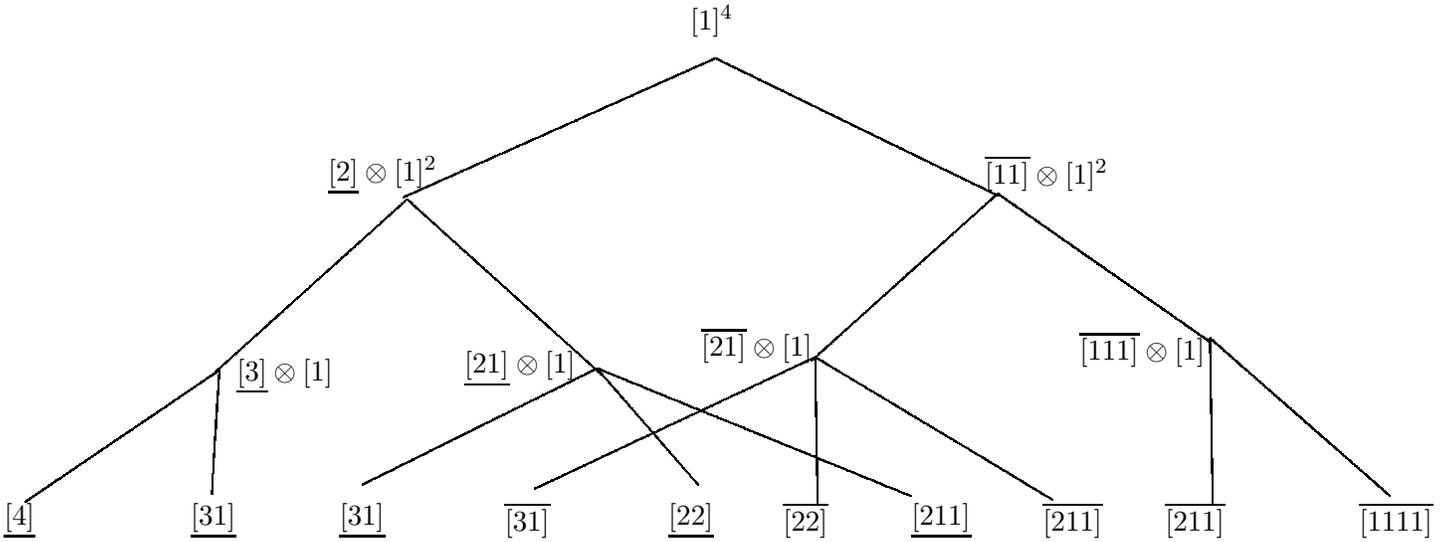
\unitlength 1mm 
\linethickness{0.4pt}
\ifx\plotpoint\undefined\newsavebox{\plotpoint}\fi 
\begin{picture}(27.5,151.5)(54,0)
\put(39,137.25){\line(1,0){.5}}
\multiput(136,141.5)(-.0755919854,-.0336976321){549}{\line(-1,0){.0755919854}}
\multiput(136,141.5)(.0693160813,-.0337338262){541}{\line(1,0){.0693160813}}
\multiput(95,122.75)(-.03739002933,-.03372434018){682}{\line(-1,0){.03739002933}}
\multiput(95,122.75)(.03731884058,-.03369565217){690}{\line(1,0){.03731884058}}
\multiput(173.5,123.5)(-.0375,-.0337121212){660}{\line(-1,0){.0375}}
\multiput(173.5,123.5)(.0482081911,-.0337030717){586}{\line(1,0){.0482081911}}
\multiput(70,100)(-.0496146435,-.0337186898){519}{\line(-1,0){.0496146435}}
\multiput(70,100)(-.0333333,-.55){30}{\line(0,-1){.55}}
\multiput(120.25,100.25)(-.0679347826,-.0336956522){460}{\line(-1,0){.0679347826}}
\multiput(149.25,101.75)(-.0717726397,-.0337186898){519}{\line(-1,0){.0717726397}}
\multiput(120.25,100.25)(.0336658354,-.0386533666){401}{\line(0,-1){.0386533666}}
\multiput(149.25,101.75)(.03125,-1.96875){10}{\line(0,-1){1.96875}}
\multiput(120.25,100.25)(.1227735369,-.05){340}{\line(1,0){.1227735369}}
\multiput(149.25,101.75)(.0831395349,-.05){380}{\line(1,0){.0831395349}}
\multiput(201.75,104.25)(.03125,-2){11}{\line(0,-1){2}}
\multiput(201.75,104.25)(.0569690265,-.05){420}{\line(1,0){.0569690265}}
\put(135.5,146.5){\makebox(0,0)[cc]{$[1]^4$}}
\put(91.75,126){\makebox(0,0)[cc]{$\underline{[2]}\otimes [1]^2$}}
\put(180,126.25){\makebox(0,0)[cc]{$\overline{[11]}\otimes [1]^2$}}
\put(78.75,99.25){\makebox(0,0)[cc]{$\underline{[3]}\otimes [1]$}}
\put(110,100.25){\makebox(0,0)[cc]{$\underline{[21]}\otimes [1]$}}
\put(141.5,103.25){\makebox(0,0)[cc]{$\overline{[21]}\otimes [1]$}}
\put(192.75,102.75){\makebox(0,0)[cc]{$\overline{[111]}\otimes [1]$}}
\put(43.5,80){\makebox(0,0)[cc]{$\underline{[4]}$}}
\put(69.25,80){\makebox(0,0)[cc]{$\underline{[31]}$}}
\put(89,80){\makebox(0,0)[cc]{$\underline{[31]}$}}
\put(111,80){\makebox(0,0)[cc]{$\overline{[31]}$}}
\put(132.75,80){\makebox(0,0)[cc]{$\underline{[22]}$}}
\put(148,80){\makebox(0,0)[cc]{$\overline{[22]}$}}
\put(166,80){\makebox(0,0)[cc]{$\underline{[211]}$}}
\put(183.5,80){\makebox(0,0)[cc]{$\overline{[211]}$}}
\put(199.75,80){\makebox(0,0)[cc]{$\overline{[211]}$}}
\put(226.5,80){\makebox(0,0)[cc]{$\overline{[1111]}$}}
\end{picture}

\vspace{-6.5cm}

At the lowest level of this picture, there emerge three representations [31] and [211],
two representations [22], while representations [4] and [1111] emerge only once. We
throw away a half of the representations that can be obtained by transposing the
corresponding Young diagrams.

Now one could construct a similar expansion, fusing at the first step, e.g.,
the two rightmost fundamental representations etc. That would lead to the same set of
representations, however, the three representations [31] emerging would be now
linear combinations of those obtained with the first pattern of fusion above.
The (orthogonal) matrix that relates representations [31] in two different
patterns of fusion is constructed from the set of the Racah coefficients
($6j$-symbols) and is exactly the $U$-matrix we are looking for.
(Similarly, we are interested in the matrix rotating the two representations [22].)

In what follows we find such matrices, constructing manifestly the vectors
in the representations and using the action of lowering and raising operators
at the tensor product:
\be\label{comult}
\downarrow\ \ \ \ \Delta(F) = q^{2H}\otimes F + F\otimes I,
\ \ \ \ \ \ \ \
\uparrow\ \ \ \ \ \Delta(E) = I \otimes E + E\otimes q^{-2H}
\ee
We begin with the highest weight $T_{++++} = T_{0000}$ and look at $1^{\otimes m}$ at
$m=2,3,4,5$.
At level one we get the pattern of the highest weights,
relevant for description of the mixing of $[m-1,1]$ representations.
Likewise, at level two we get the description of mixing in the $[m-2,2]$ sector.

\newpage

\subsection{$m=2$}

In this case, there is only one pattern of fusion, the vectors in the representations
being manifestly indicated at the picture:

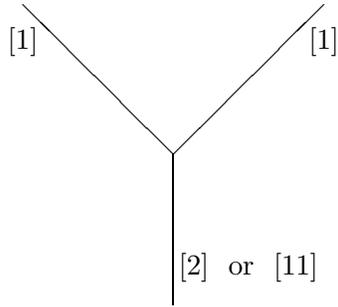
\vspace{1cm}

\begin{picture}(20,10)(-30,10)
\put(0,0){\line(0,-1){20}}
\put(0,0){\line(-1,1){20}}
\put(0,0){\line(1,1){20}}
\put(-20,15){\makebox(0,0)[cc]{$[1]$}}
\put(20,15){\makebox(0,0)[cc]{$[1]$}}
\put(10,-15){\makebox(0,0)[cc]{$[2]$ \ or \ $[11]$}}
\end{picture}

\be
\hspace{8cm}\begin{array}{ccc}
&[1]^2& \\
\swarrow\!\!\!\!\!\!\!\!\!\!\!\!\!\!\!\!\!\!\!\!\!\!\!\!
&&\!\!\!\!\!\!\!\!\!\!\!\!\!\!\!\!\!\!\!\!\!\!\!\! \searrow \\
\underline{[2]} && \overline{[11]} \\
T_{00}&&\\
T_{10} + qT_{01} &&  -qT_{10} + T_{01}\\
T_{11} &&
\end{array}
\ee

\bigskip

\subsection{$m=3$}

In this case, there are two possible fusion trees:

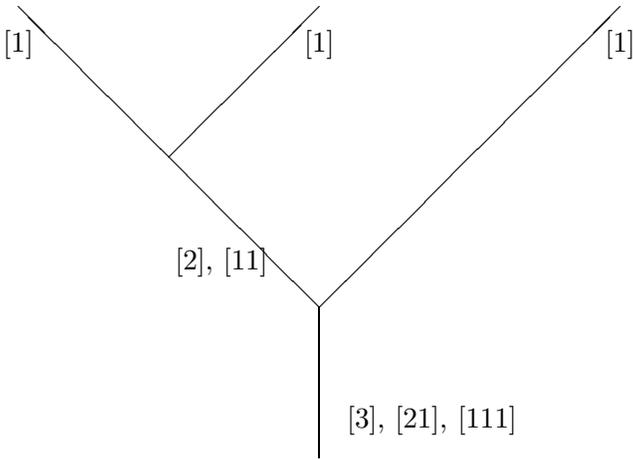
\begin{picture}(10,20)(-30,30)
\put(0,0){\line(0,-1){20}}
\put(0,0){\line(-1,1){40}}
\put(-20,20){\line(1,1){20}}
\put(0,0){\line(1,1){40}}
\put(-40,35){\makebox(0,0)[cc]{$[1]$}}
\put(40,35){\makebox(0,0)[cc]{$[1]$}}
\put(0,35){\makebox(0,0)[cc]{$[1]$}}
\put(-13,6){\makebox(0,0)[cc]{$[2]$, $[11]$}}
\put(15,-15){\makebox(0,0)[cc]{$[3]$, $[21]$, $[111]$}}
\end{picture}

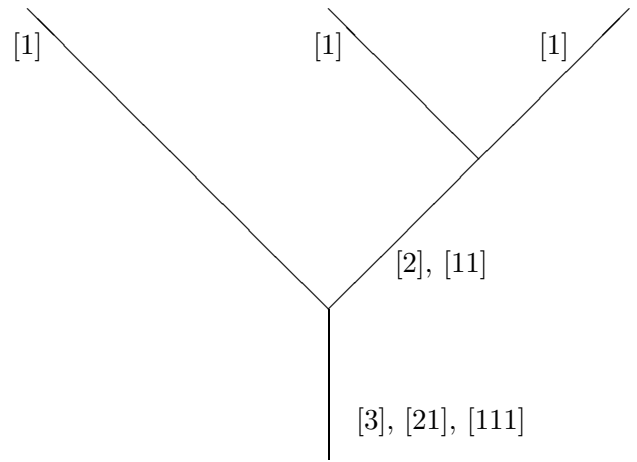
\begin{picture}(20,-30)(-130,25.5)
\put(0,0){\line(0,-1){20}}
\put(0,0){\line(-1,1){40}}
\put(20,20){\line(-1,1){20}}
\put(0,0){\line(1,1){40}}
\put(-40,35){\makebox(0,0)[cc]{$[1]$}}
\put(30,35){\makebox(0,0)[cc]{$[1]$}}
\put(0,35){\makebox(0,0)[cc]{$[1]$}}
\put(15,6){\makebox(0,0)[cc]{$[2]$, $[11]$}}
\put(15,-15){\makebox(0,0)[cc]{$[3]$, $[21]$, $[111]$}}
\end{picture}

\vspace{6cm}

\noindent
Let us list the vectors in the representations emerging in the right tree through fusions
in the $SU(2)$ case (so that [111] does not emerge):

\be
{\footnotesize
\begin{array}{cccc}
&\hspace{2cm}\l[1]^3 &&\\ \\&\hspace{-3cm}\swarrow&\searrow&\\ \\
\l[1]\otimes\underline{[2]}& && [1]\otimes\overline{[11]} \\
T_{000}  &&& \\
T_{001};\ \ \ \ \ \ T_{100}+qT_{010} &
&& -qT_{100} + T_{010} \\
T_{110};\ \ \ \ \ \ T_{101} + qT_{011} &&& -qT_{101}+T_{011} \\
T_{111} &&& \\
\\  \downarrow&\hspace{-1cm}\searrow &&\downarrow\\
\\
\underline{[3]} & \underline{[21]} && \overline{[21]} \\
T_{000} &&&\\
T_{100}+qT_{010} + q^2T_{001} & qT_{100} + q^2T_{010} - [2]_qT_{001} &
\hspace{3cm}& -qT_{100} + T_{010}\\
T_{110}+qT_{101} + q^2T_{011} & q^2[2]_qT_{110} - T_{101}-qT_{011} &&
-qT_{101} + T_{011} \\
T_{111} &&&
\end{array}
}\nn
\ee

We did not normalize the vectors here.

One can see at this scheme how the vectors at each grade are parted between representations.
For instance, at the first level of the hierarchy one just has one grade zero vector
$T_{000}$, three grade one vectors $T_{001}$, $T_{010}$ and $T_{100}$, three grade two
vectors $T_{011}$, $T_{110}$ and $T_{101}$ and one grade three vector $T_{111}$. These
vectors are further distributed between various representations. At the second level
of the hierarchy, $T_{000}$ goes to $[1]\otimes [2]$, while
the three grade one vectors, $T_{001}$, $T_{010}$ and $T_{100}$ are distributed between
$[1]\otimes [2]$ and $[1]\otimes [11]$ etc.

From now on, we are interested in the vectors in irreducible representations only, i.e.
at the last level of the hierarchy. In particular, the content of irreps for the left tree is
\be
[1]^3\to \underline{[2]}\otimes [1]\to \underline{[3]}:\hspace{3cm}
\begin{array}{c}
T_{000}\\
T_{100}+qT_{010} + q^2T_{001}\\
T_{110}+qT_{101} + q^2T_{011}\\
T_{111}
\end{array}\nn
\ee

\be
[1]^3\to \underline{[2]}\otimes [1]\to \underline{[21]}:\hspace{3cm}
\begin{array}{c}
-q^2[2]_qT_{100} + T_{010} +qT_{001}\\
-[2]_qT_{011}+q^2T_{101}+qT_{110}
\end{array}\nn
\ee

\be
[1]^3\to \overline{[11]}\otimes [1]\to \overline{[21]}:\hspace{3cm}
\begin{array}{c}
-qT_{010} + T_{001}\\
-qT_{110} + T_{101}
\end{array}\nn
\ee

Looking at the elements of the representations at the first grade,
one immediately read out expressions for the highest weights
of the two representations $[21]$ at the left and right trees accordingly:
\be
\hbox{left\ tree}\longrightarrow\hspace{0.2cm}\left\{\begin{array}{ccl}
\underline{[21]}:& -q^2[2]_qT_{100} + T_{010} +qT_{001};& \hspace{1cm}||\ldots||^2=
q^3[2]_q[3]_q\\
\\
\overline{[21]}:& -qT_{010} + T_{001};&\hspace{1cm}||\ldots||^2= q[2]_q
\end{array}\right.\nn
\ee
\be
\hbox{right\ tree}\longrightarrow\hspace{0.2cm}\left\{\begin{array}{ccl}
\underline{[21]}:& qT_{100} + q^2T_{010} -[2]_qT_{001};&\hspace{1cm}||\ldots||^2=
q[2]_q[3]_q\\
\\
\overline{[21]}: &-qT_{100} + T_{010};&\hspace{1cm} ||\ldots||^2=q[2]_q
\end{array}\right.\nn
\ee
The last column contains the squares of norms of these highest weights
(the sums of squares of the coefficients\footnote{
Throughout the paper we
assume that the norms are defined in the naive way: :
$$
\left< T_{k_1k_2\ldots } \Big|
T_{k_1'k_2'\ldots}\right> = \prod_{i=1}^m \delta_{k_ik_i'}
$$
In fact, the scalar product is a little more complicated,
if generators $E$ and $F$ with comultiplication (\ref{comult})
are conjugate (the representation is unitary). Indeed,
$$
\left< \Delta(E) (A\otimes B) | C\otimes D \right>
= \left< A\otimes B | \Delta(F) (C\otimes D)\right>
$$
implies that
$$
\left<A\otimes B | C\otimes D\right> = q^{2h_Ah_D}<A|C><B|D>
$$
and in general
$$
\left<\otimes_{i=1}^m A_i \Big| \otimes_{i=1}^m C_i\right>
= q^{2\sum_{i<j}h_{A_i}h_{C_j}} \prod_{i=1}^m <A_i|C_i>
$$
However, for $R=[1]$ the "cocycle" factor,
though non-trivial, can be neglected in all our calculations,
because in this case of $p=1$ (only!)
the norms of all the elements of the same grade
in Verma module are the same:
$||T_{k_1k_2\ldots}||$ depends only on the sum $\sum_i k_i$
and not on a particular choice of $\{k_i\}$.
}).
If highest weights are divided by these norms,
the mixing matrices are orthogonal.
This is the way we define them in what follows. In particular, in this way
one immediately obtains the orthogonal $2\times 2$ mixing matrix $U$:
\be
U = \ldb\begin{array}{cc} {1\over [2]_q} & {\sqrt{[3]_q}\over [2]_q} \\
-{\sqrt{[3]_q}\over [2]_q} & {1\over [2]_q}\end{array}\rdb
\ee
which is exactly (\ref{2sU}).

\subsection{$m=4$}

In this case, there are four possible fusion trees.

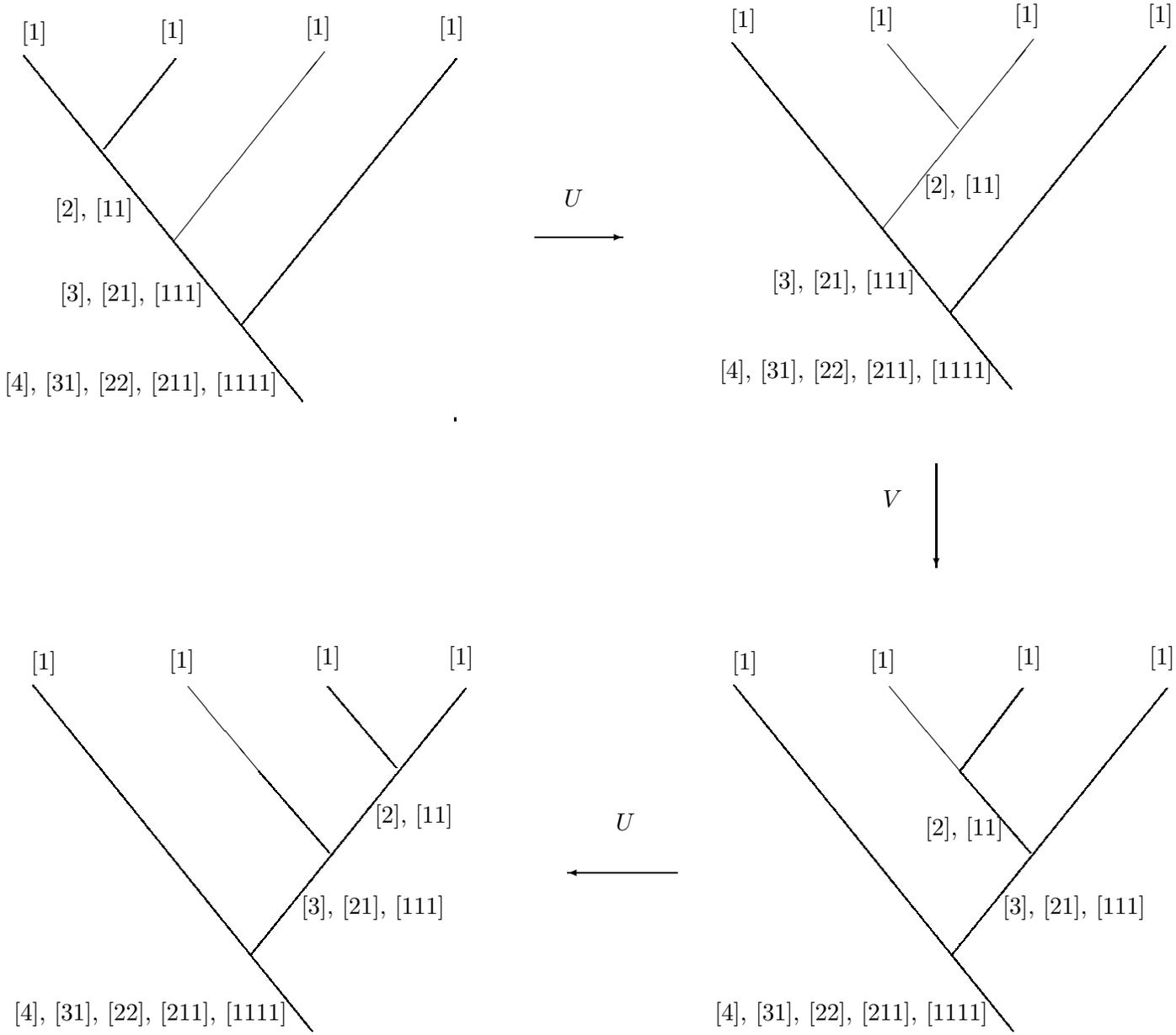
\unitlength 0.99mm 
\linethickness{0.4pt}
\ifx\plotpoint\undefined\newsavebox{\plotpoint}\fi 
\begin{picture}(197.75,233)(145,-150)
\put(204.25,18.5){\framebox(0,.25)[]{}}
\put(159.5,47.25){\line(4,5){24}}
\multiput(159.75,76.25)(-.0337243402,-.0425219941){341}{\line(0,-1){.0425219941}}
\multiput(135.5,76.75)(.03373768006,-.04188779378){1319}{\line(0,-1){.04188779378}}
\put(171.5,35.5){\line(0,-1){.25}}
\multiput(204.5,76.25)(-.03371062992,-.04183070866){1016}{\line(0,-1){.04183070866}}
\put(137.25,80.25){\makebox(0,0)[cc]{$[1]$}}
\put(159.25,80.5){\makebox(0,0)[cc]{$[1]$}}
\put(182.5,81){\makebox(0,0)[cc]{$[1]$}}
\put(203.75,81){\makebox(0,0)[cc]{$[1]$}}
\put(146.75,52){\makebox(0,0)[cc]{$[2]$, $[11]$}}
\put(152.75,38.5){\makebox(0,0)[cc]{$[3]$, $[21]$, $[111]$}}
\put(154.25,24.25){\makebox(0,0)[cc]{$[4]$, $[31]$, $[22]$, $[211]$, $[1111]$}}
\put(272.5,49.25){\line(4,5){24}}
\multiput(248.5,78.75)(.03373768006,-.04188779378){1319}{\line(0,-1){.04188779378}}
\put(284.5,37.5){\line(0,-1){.25}}
\multiput(317.5,78.25)(-.03371062992,-.04183070866){1016}{\line(0,-1){.04183070866}}
\put(250.25,82.25){\makebox(0,0)[cc]{$[1]$}}
\put(272.25,82.5){\makebox(0,0)[cc]{$[1]$}}
\put(295.5,83){\makebox(0,0)[cc]{$[1]$}}
\put(316.75,83){\makebox(0,0)[cc]{$[1]$}}
\put(266.25,40.5){\makebox(0,0)[cc]{$[3]$, $[21]$, $[111]$}}
\put(268.25,26.25){\makebox(0,0)[cc]{$[4]$, $[31]$, $[22]$, $[211]$, $[1111]$}}
\put(273.25,78.5){\line(5,-6){11.25}}
\put(285.25,55.75){\makebox(0,0)[cc]{$[2]$, $[11]$}}
\put(217,47.75){\vector(1,0){14}}
\put(223.25,54){\makebox(0,0)[cc]{${\huge U}$}}
\multiput(248.75,-23.75)(.03373768006,-.04188779378){1319}{\line(0,-1){.04188779378}}
\multiput(137,-23.75)(.03373768006,-.04188779378){1319}{\line(0,-1){.04188779378}}
\put(284.75,-65){\line(0,-1){.25}}
\put(173,-65){\line(0,-1){.25}}
\multiput(317.75,-24.25)(-.03371062992,-.04183070866){1016}{\line(0,-1){.04183070866}}
\multiput(206,-24.25)(-.03371062992,-.04183070866){1016}{\line(0,-1){.04183070866}}
\put(250.5,-20.25){\makebox(0,0)[cc]{$[1]$}}
\put(138.75,-20.25){\makebox(0,0)[cc]{$[1]$}}
\put(272.5,-20){\makebox(0,0)[cc]{$[1]$}}
\put(160.75,-20){\makebox(0,0)[cc]{$[1]$}}
\put(295.75,-19.5){\makebox(0,0)[cc]{$[1]$}}
\put(184,-19.5){\makebox(0,0)[cc]{$[1]$}}
\put(317,-19.5){\makebox(0,0)[cc]{$[1]$}}
\put(205.25,-19.5){\makebox(0,0)[cc]{$[1]$}}
\put(285.5,-46.75){\makebox(0,0)[cc]{$[2]$, $[11]$}}
\put(197.75,-44.75){\makebox(0,0)[cc]{$[2]$, $[11]$}}
\put(303,-59.25){\makebox(0,0)[cc]{$[3]$, $[21]$, $[111]$}}
\put(191.25,-59.25){\makebox(0,0)[cc]{$[3]$, $[21]$, $[111]$}}
\put(267.5,-76.25){\makebox(0,0)[cc]{$[4]$, $[31]$, $[22]$, $[211]$, $[1111]$}}
\put(155.75,-76.25){\makebox(0,0)[cc]{$[4]$, $[31]$, $[22]$, $[211]$, $[1111]$}}
\put(273.5,-24){\line(5,-6){11.25}}
\put(161.75,-24){\line(5,-6){11.25}}
\put(281,11.5){\vector(0,-1){16.5}}
\put(274,6){\makebox(0,0)[cc]{$V$}}
\multiput(294.75,-24.25)(-.0336700337,-.0446127946){297}{\line(0,-1){.0446127946}}
\multiput(284.75,-37.5)(.0336826347,-.0389221557){334}{\line(0,-1){.0389221557}}
\multiput(173,-37.5)(.0336826347,-.0389221557){334}{\line(0,-1){.0389221557}}
\multiput(184,-24)(.0336391437,-.0397553517){327}{\line(0,-1){.0397553517}}
\put(239.75,-53.75){\vector(-1,0){17.5}}
\put(231.5,-45.5){\makebox(0,0)[cc]{$U$}}
\end{picture}

\vspace{-6cm}

Since we deal here with $SU(2)$ only, representations $[111], [211]$ and
$[1111]$ are not seen, however, their transposed are. Hence, the complete results
in the $4$-strand case can be obtained from the $SU(2)$ group only (and the fundamental
representations). The vectors in the representations in the first tree are

\be
[1]^4\to \underline{[2]}\otimes [1]^2\to \underline{[3]}\times[1]\to
\underline{[4]}:\hspace{1.5cm}
\begin{array}{c}
T_{0000}\\
T_{1000}+qT_{0100}+q^2T_{0010}+q^3T_{0001}\\
T_{1100}+qT_{1010}+q^2T_{0110}
+q^2T_{1001} + q^3T_{0101} + q^4T_{0011}\\
T_{1110} + qT_{1101} + q^2T_{1011} + q^3T_{0111}\\
T_{1111}
\end{array}\nn
\ee

\be
[1]^4\to \underline{[2]}\otimes [1]^2\to \underline{[3]}\times[1]\to
\underline{[31]}:\hspace{1.5cm}
\begin{array}{c}
T_{1000} +qT_{0100} + q^2T_{0010} -\frac{[3]_q}{q}T_{0001}\\
T_{1100}+qT_{1010}+q^2T_{0110} -
-\frac{1}{q^2}T_{1001} -\frac{1}{q}T_{0101} - T_{0011}\\
T_{1110} -\frac{1}{q^3}T_{1101} - \frac{1}{q^2}T_{1011} - \frac{1}{q}T_{0111}
\end{array}\nn
\ee

\be
[1]^4\to \underline{[2]}\otimes [1]^2\to \underline{[21]}\times[1]\to
\underline{[31]}:\hspace{1cm}
\begin{array}{c}
T_{1000} + qT_{0100} - \frac{[2]_q}{q}T_{0010}\\
\l[2]_qT_{1100} - \frac{1}{q^2}T_{1010}-\frac{1}{q}T_{0110}
+qT_{1001}  + q^2T_{0101} -[2]_qT_{0011}\\
\l[2]_qT_{1101}-\frac{1}{q^2}T_{1011}-\frac{1}{q}T_{0111}
\end{array}\nn
\ee

\be
[1]^4\to \overline{[11]}\otimes [1]^2\to \overline{[21]}\times[1]\to
\overline{[31]}:\hspace{1.5cm}
\begin{array}{c}
-qT_{1000}+T_{0100}\\
-qT_{1010} + T_{0110}
-q^2T_{1001} + qT_{0101}\\
-qT_{1011} + T_{0111}
\end{array}\nn
\ee

Looking at the elements of the representations at the first grade,
one immediately reads off expressions for the three highest weights
of the three representations $[31]$:
\be
\begin{array}{ccc}
\underline{[31]} &  T_{1000} +qT_{0100} + q^2T_{0010} -\frac{[3]_q}{q}T_{0001}
\ \ \ \ \ \ \
& q^{-1}[2]_q[3]_q(q^2+q^{-2})  \\
\underline{[31]} & T_{1000} + qT_{0100} - \frac{[2]_q}{q}T_{0010} & q^{-1}[2]_q[3]_q \\
\overline{[31]} &  -qT_{1000}+T_{0100} & q[2]_q
\end{array}\nn
\ee

At the second grade the four representations $[4] + 3\,[31]$ do not exhaust all
of the tensor $T$, even for $SU(2)$.
There are actually six elements of $T$ at this grade, and only four belong
to those four representations.
Thus, there are two more highest weights: those of representations $[22]$
(as to $[211]$ and $[1111]$, they are "non-observable" while we deal with $SU(2)$).
These two new weights are:
\be
\begin{array}{ccccc}
\underline{[21]}\otimes [1] & \longrightarrow & \underline{[22]}
& \ \ \ \
q[2]_qT_{1100} - \frac{1}{q}T_{1010}-T_{0110}
-  T_{1001}  - qT_{0101} +\frac{[2]_q}{q}T_{0011}
\ \ \ \ \ \ \ & [2]_q^2[3]_q \\ \\
\overline{[21]}\otimes [1] & \longrightarrow & \overline{[22]}
&  -qT_{1010} + T_{0110} + T_{1001} - \frac{1}{q}T_{0101}   &  [2]_q^2
\end{array}\nn
\ee

Now we switch to the other tree related with the just considered one by the
mixing matrix $U$ (see the trees in the diagram above) and obtain the following
highest vectors
in the representations $[31]$:

\be
\begin{array}{ccc}
\underline{[31]} &  T_{1000} +qT_{0100} + q^2T_{0010} -\frac{[3]_q}{q}T_{0001}
& q^{-1}[2]_q[3]_q(q^2+q^{-2})  \\
\underline{[31]} & -q^2[2]_qT_{1000} + T_{0100} + qT_{0010}  & q^3[2]_q[3]_q \\
\overline{[31]} &  -qT_{0100}+T_{0010} & q[2]_q
\end{array}\nn
\ee

Applying further the matrix $V$, we come to the next tree with the highest vectors
in the representations $[31]$ as follows:

\be
\begin{array}{ccc}
\underline{[31]} & -q^2[3]_q T_{1000} + T_{0100} + qT_{0010} +q^2T_{0001} & q^5[2]_q[3]_q(q^2+q^{-2}) \\
\underline{[31]} &  T_{0100} + qT_{0010} -\frac{[2]_q}{q}T_{0001} & q^{-1}[2]_q[3]_q  \\
\overline{[31]} &  -qT_{0100}+T_{0010} & q[2]_q
\end{array}\nn
\ee

The last tree can be obtained with action of the same matrix $U$, the highest vectors
being

\be
\begin{array}{ccc}
\underline{[31]} & -q^3[3]_q T_{1000} + T_{0100} + qT_{0010} +q^2T_{0001} & q^5[2]_q[3]_q(q^2+q^{-2})   \\
\underline{[31]} & -q^2[2]_q T_{0100} + T_{0010} + qT_{0001}  & q^{-1}[2]_q[3]_q  \\
\overline{[31]} &  -qT_{0010} + T_{0001} & q[2]_q
\end{array}\nn
\ee

Similarly to the previous subsection, from these manifest expressions for the vectors
one can read off the mixing matrices of rotation:
\be
U = \ldb\begin{array}{ccc}
1 & 0 & 0 \\
0 & \frac{1}{[2]_q} & \frac{\sqrt{[3]_q}}{[2]_q} \\
0 & -\frac{\sqrt{[3]_q}}{[2]_q} & \frac{1}{[2]_q}
\end{array}\rdb
\ee
\be
V = \ldb\begin{array}{ccc}
 \frac{1}{[3]_q} & \frac{[2]_q\sqrt{q^2+q^{-2}}}{[3]_q} & 0 \\
 -\frac{[2]_q\sqrt{q^2+q^{-2}}}{[3]_q} & \frac{1}{[3]_q} & 0  \\
0 & 0 & 1
\end{array}\rdb
\ee
so that
\be\label{R31}
{\mathfrak{R}}^{[31]}_1 = \ldb\begin{array}{ccc}
q & & \\ & q & \\ &&-\frac{1}{q}
\end{array}\rdb,   \\
{\mathfrak{R}}_2^{[31]} = U\ {\mathfrak{R}}_1^{[31]}\ \tilde U  =
\ldb\begin{array}{ccc}
q &0 &0 \\0 & -\frac{1}{q^2[2]_q} & -\frac{\sqrt{[3]_q}}{[2]_q} \\
0&-\frac{\sqrt{[3]_q}}{[2]_q} &\frac{q^2}{[2]_q}
\end{array}\rdb,   \\
{\mathfrak{R}}_3^{[31]} = UV\ {\mathfrak{R}}_2^{[31]}\ \widetilde{UV}
= \ldb\begin{array}{ccc}
-\frac{1}{q^3[3]_q} & -\frac{[2]_q\sqrt{q^2+q^{-2}}}{[3]_q} &0 \\
-\frac{[2]_q\sqrt{q^2+q^{-2}}}{[3]_q} & \frac{q^3}{[3]_q} & 0 \\ 0&0&q
\end{array}\rdb
\ee
and
\be\label{R22}
{\mathfrak{R}}_1^{[22]} = \ldb\begin{array}{cc} q & \\ & -\frac{1}{q} \end{array}\rdb,   \\
{\mathfrak{R}}_2^{[22]} = U^{[22]}\ {\mathfrak{R}}_1^{[22]}\ \tilde U^{[22]} =
\ldb\begin{array}{cc} -\frac{1}{q^2[2]_q} & -\frac{\sqrt{[3]_q}}{[2]_q} \\
-\frac{\sqrt{[3]}}{[2]_q} & \frac{q^2}{[2]_q} \end{array}\rdb, \\
{\mathfrak{R}}_3^{[22]} = \tilde U^{[22]}\ {\mathfrak{R}}_2^{[22]} \ U^{[22]} =
\ldb\begin{array}{cc} q & \\ & -\frac{1}{q} \end{array}\rdb,\ \hbox{since}
\hspace{.5cm} V^{[22]}=\left(U^{[22]}\right)^{-2}=
\ldb\begin{array}{cc} -\frac{[4]}{[2]^3}&-2\frac{\sqrt{[3]}}{[2]^2}\\
2\frac{\sqrt{[3]}}{[2]^2}&-\frac{[4]}{[2]^3} \end{array}\rdb
\ee
where tilde implies the transposed matrix (or inverse, since all they are
orthogonal).

\subsection{$m=5$}

For $m=5$ there is one transition between the trees of a new type, $W$:

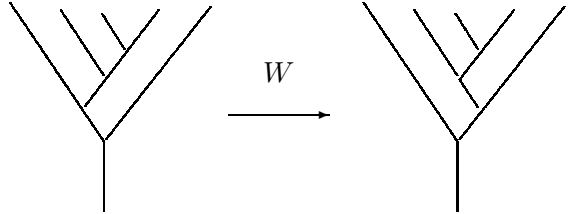
\unitlength 1mm 
\linethickness{0.4pt}
\ifx\plotpoint\undefined\newsavebox{\plotpoint}\fi 
\begin{picture}(121.25,44.75)(0,160)
\multiput(47.5,194.75)(.0336927224,-.0498652291){371}{\line(0,-1){.0498652291}}
\multiput(94.5,194.75)(.0336927224,-.0498652291){371}{\line(0,-1){.0498652291}}
\put(60,176.25){\line(0,-1){9.25}}
\put(107,176.25){\line(0,-1){9.25}}
\multiput(60.25,176.5)(.0337349398,.0427710843){415}{\line(0,1){.0427710843}}
\multiput(107.25,176.5)(.0337349398,.0427710843){415}{\line(0,1){.0427710843}}
\multiput(67.25,193.75)(-.0337370242,-.0441176471){289}{\line(0,-1){.0441176471}}
\multiput(59.75,193.25)(.03370787,-.05337079){89}{\line(0,-1){.05337079}}
\multiput(106.75,193.25)(.03370787,-.05337079){89}{\line(0,-1){.05337079}}
\multiput(54.25,193.75)(.033625731,-.051169591){171}{\line(0,-1){.051169591}}
\multiput(101.25,193.75)(.033625731,-.051169591){171}{\line(0,-1){.051169591}}
\put(76.75,179.75){\line(-1,0){.25}}
\put(76.75,179.75){\vector(1,0){13.25}}
\put(83.25,185.5){\makebox(0,0)[cc]{$W$}}
\multiput(114.5,193.75)(-.03372093,-.043023256){215}{\line(0,-1){.043023256}}
\multiput(107.25,184.5)(.03333333,-.05){75}{\line(0,-1){.05}}
\end{picture}

The highest vectors in the four representations [41] in these two trees
correspondingly look like

\be
\begin{array}{ccc}
\underline{[41]} & T_{10000} + qT_{01000} + q^2T_{00100} + q^3T_{00010} -
\frac{[4]_q}{q}T_{00001}
& q^{-1}[2]_q[5]_q(q^2+q^{-2})   \\
\underline{[41]} & -q^3[3]_qT_{10000} + T_{01000} + qT_{00100} + q^2T_{00010} &
q^5[2]_q[3]_q(q^2+q^{-2})   \\
\underline{[41]} & -q^2[2]_qT_{01000}+T_{00100} + qT_{00010} & q^3[2]_q[3]_q   \\
\overline{[41]} & -qT_{00100} + T_{00010} & q[2]_q
\end{array}\nn
\ee

and

\be
\begin{array}{ccc}
\underline{[41]} & -q^4[4]_qT_{10000} + T_{01000} + qT_{00100} + q^2T_{00010} +
q^3T_{00001}
& q^{7}[2]_q[5]_q(q^2+q^{-2})   \\
\underline{[41]} & T_{01000} + qT_{00100} + q^2T_{00010} - \frac{[3]_q}{q}T_{00001}
& q^{-1}[2]_q[3]_q(q^2+q^{-2})   \\
\underline{[41]} &-q^2[2]_qT_{01000}+T_{00100} + qT_{00010} & q^3[2]_q[3]_q   \\
\overline{[41]} & -qT_{00100} + T_{00010} & q[2]_q
\end{array}\nn
\ee

They are related by the new $4\times 4$ matrix $W$:
\be
W = \ldb\begin{array}{cccc}
 \frac{1}{[4]_q} & \frac{\sqrt{[3]_q[5]_q}}{[4]_q} & 0 & 0 \\
 -\frac{\sqrt{[3]_q[5]_q}}{[4]} & \frac{1}{[4]_q} & 0 & 0 \\
0 & 0 & 1 & 0 \\
0 & 0 & 0 & 1
\end{array}\rdb
\ee

In general in the sector $[m-1,1] \in [1]^{\otimes m}$
one has $(m-1)\times (m-1)$ matrices
\be
U_{(k-1)} = \ldb\begin{array}{cccccccccc}
1 & 0 & \ldots & 0 & 0&0 & 0 & \ldots & 0 & 0\\
0 & 1 & &0&0&0&0& 0 & 0 \\
\ldots &&&&&&&&& \\
0 & 0 && 1 & 0&0 & 0 && 0 & 0 \\
0& 0 && 0  & \displaystyle{ \frac{1}{[k]_q}} & \displaystyle{{\sqrt{[k]_q^2-1}\over [k]_q}}
 & 0 && 0 & 0\\
0 & 0 && 0 & \displaystyle{-{\sqrt{[k]_q^2-1}\over [k]_q}}&
\displaystyle{\frac{1}{[k]_q}} & 0 && 0 & 0 \\
0 & 0 && 0 & 0&0 & 1 && 0 & 0 \\
\ldots &&&&&&&&& \\
0 & 0 && 0 & 0&0 & 0 && 1 & 0 \\
0 & 0 && 0 & 0&0 & 1 && 0 & 1
\end{array}\rdb
\ee
where non-unit $2\times 2$ matrix
stands in the rows and columns $m-k, m-k+1$ of the full matrix,
$k = 2,3,\ldots,m-1$.

In the previous examples $U = U_{(1)}$, $V = U_{(2)}$, $W= U_{(3)}$:
\be\label{UV41}
U = \ldb\begin{array}{cccccc}
1 && 0 & 0 & 0 & 0 \\
& \ldots &  &  &  &  \\
0 && 1 & 0 & 0 & 0\\
0& & 0 & 1 & 0 & 0 \\
0 && 0 & 0 & \frac{1}{[2]_q} & \frac{\sqrt{[3]_q}}{[2]_q} \\
0 && 0 & 0 & -\frac{\sqrt{[3]_q}}{[2]_q} & \frac{1}{[2]_q}
\end{array}\rdb,
\\
V = \ldb\begin{array}{cccccc}
1 && 0 & 0 & 0 & 0 \\
& \ldots &  &  &  &  \\
0 && 1 & 0 & 0 & 0\\
0 && 0 & \frac{1}{[3]_q} & \frac{[2]_q\sqrt{q^2+q^{-2}}}{[3]_q} & 0 \\
0 && 0 & -\frac{[2]_q\sqrt{q^2+q^{-2}}}{[3]_q} & \frac{1}{[3]_q} & 0  \\
0 && 0 & 0 & 0 & 1
\end{array}\rdb
\ee
\be\label{W41}
W = \ldb\begin{array}{cccccc}
1 && 0 & 0 & 0 & 0 \\
& \ldots &  &  &  &  \\
0 & 0 & \frac{1}{[4]_q} & \frac{\sqrt{[3]_q[5]_q}}{[4]_q} & 0 & 0 \\
0 & 0 & -\frac{\sqrt{[3]_q[5]_q}}{[4]_q} & \frac{1}{[4]_q} & 0 & 0 \\
0 & 0 & 0 & 0 & 1 & 0 \\
0 & 0 & 0 & 0 & 0 & 1
\end{array}\rdb
\ee
since
\be
\l[2]_q^2-1 = [3]_q,   \\
\l[3]_q^2 - 1 = [2]_q^2[2]_{q^2}=[2]_q[4]_q,   \\
\l[4]_q^2-1 = [3]_q[5]_q,   \\
\l[5]_q^2-1 = [2]_q[3]_q[2]_{q^2}[2]_{q^3}=[4]_q[6]_q, \\
\ldots,\\
\l[k]_q^2-1=[k+1]_q[k-1]_q
\ee

\section{4-strand braids}

In order to calculate $H_{[1]}$ in the case of four strands, one needs to deal with traces
\be
\tr_{R} \hat{\mathfrak{R}}_{1}^{a_1}\hat{\mathfrak{R}}_{2}^{b_1}\hat{\mathfrak{R}}_{3}^{c_1}
\hat{\mathfrak{R}}_{1}^{a_2}\hat{\mathfrak{R}}_{2}^{b_2}\hat{\mathfrak{R}}_{3}^{c_2}\ldots
\ee
over five possible
representations $R$ emerging in the decomposition
\be
[1]^4=[4]+3[31]+2[22]+3[211]+[1111]
\ee
i.e. with the sum
\be
H_{[1]}=
\Tr \left\{\Big({\mathfrak{R}}\otimes I\otimes I\Big)^{a_1}
\Big(I \otimes {\mathfrak{R}}\otimes I\Big)^{b_1}
\Big(I \otimes I\otimes {\mathfrak{R}}\Big)^{c_1}
\Big({\mathfrak{R}}\otimes I\otimes I\Big)^{a_2} \Big(I \otimes {\mathfrak{R}}
\otimes I\Big)^{b_2}
\Big(I \otimes I\otimes {\mathfrak{R}}\Big)^{c_2} \ldots\right\}=
\\
=\sum_{{Q=[1111],[211],}\atop{[22],[31],[4]}}\tr\left\{
\left(\hat{\mathfrak{R}}_1^Q\right)^{a_1}
\left(\hat{\mathfrak{R}}_2^Q\right)^{b_1}
\left(\hat{\mathfrak{R}}_3^Q\right)^{c_1}\left(\hat{\mathfrak{R}}_1^Q\right)^{a_2}
\left(\hat{\mathfrak{R}}_2^Q\right)^{b_2}
\left(\hat{\mathfrak{R}}_3^Q\right)^{c_2}\ldots\right\}=\\
=\sum_{{Q=[1111],[211],}\atop{[22],[31],[4]}}\tr\left\{
\left(\hat{\mathfrak{R}}_1^Q\right)^{a_1}{\hat U^Q}
\left(\hat{\mathfrak{R}}_1^Q\right)^{b_1}{\hat V^Q}{\hat U^Q}
\left(\hat{\mathfrak{R}}_1^Q\right)^{c_1}\widetilde{{\hat U^Q}{\hat V^Q}{\hat U^Q}}
\times\right.\\\left.\times
\left(\hat{\mathfrak{R}}_1^Q\right)^{a_2}{\hat U^Q}
\left(\hat{\mathfrak{R}}_1^Q\right)^{b_2}{\hat V^Q}{\hat U^Q}
\left(\hat{\mathfrak{R}}_1^Q\right)^{c_2}\widetilde{{\hat U^Q}{\hat V^Q}{\hat U^Q}}
\ldots\right\}
\ee
since
\be
\hat{\mathfrak{R}}_2^Q={\hat U^Q}\hat{\mathfrak{R}}_1^Q\widetilde{{\hat U^Q}}\ ,
\hspace{2cm}
\hat{\mathfrak{R}}_3^Q={\hat U^Q}{\hat V^Q}{\hat U^Q}
\hat{\mathfrak{R}}_1^Q\widetilde{{\hat U^Q}{\hat V^Q}{\hat U^Q}}
\ee
The corresponding manifest expressions for ${\mathfrak{R}}$- and $U$-, $V$-matrices
can be found in the previous section.

In particular, the contributions from $[4]$ and $[1111]$ sectors to the HOMFLY polynomials are trivial:
there are no mixing matrices and all what we get is simply
\be
q^{a_1+b_1+c_1+a_2+b_2+c_2+\ldots}\ S_4^* +
(-1/q)^{a_1+b_1+c_1+a_2+b_2+c_2+\ldots}\ S_{1111}^*
\ee
Evaluation of contribution from $[22]$ is straightforward:
the $2\times 2$ ${\mathfrak{R}}$-matrices are actually given by the same mixing matrix
as that in the $[21]$ sector
for the $3$-strand braids, hence, by the same ${\mathfrak{R}}$-matrix at
the adjacent strand, (\ref{R22}b)
and trivial ${\mathfrak{R}}$-matrix at the non-adjacent strand (\ref{R22}c).

However, the contribution to, say, torus knots is different,
because now the $\hat{\mathfrak{R}}$-matrix enters three times instead
of two in the $[21]$ case, and
\be
\tr_{[22]} \Big(\hat{\mathfrak{R}}_1\hat{\mathfrak{R}}_2\hat{\mathfrak{R}}_3\Big)^n
= 1^n + (-1)^n = \left\{\begin{array}{ccc}
0 & {\rm for \ odd} & n \neq 0\ {\rm mod} \ 2 \\
2 & {\rm for \ even} & n = 0\ {\rm mod} \ 2
\end{array}\right.
\ee

Similarly, the $[31]$-contribution is done with using the $\hat{\mathfrak{R}}$-matrices from (\ref{R31}).
For torus knots
\be
\tr_{[31]} \Big(\hat{\mathfrak{R}}_1\hat{\mathfrak{R}}_2\hat{\mathfrak{R}}_3\Big)^n
=\left.\sum_{k=1..L-1}e^{2\pi ikn\over L}q^n\right|_{L=4}
=\left\{ \begin{array}{ccc}
-q^n & {\rm for} & n\neq 0\ {\rm mod}\ 4 \\
3q^n & {\rm for} & n=0\ {\rm mod}\ 4
\end{array}\right.
\ee

The contribution of the $[211]$
sector can be obtained from that of $[31]$ by a simple substitution
$q \rightarrow (-1/q)$ (in particular, $S_{211}^*(q) = S_{31}^*(-1/q)$).

\newpage

\thispagestyle{empty}

\hspace{-1cm}
{\tiny{\rotate{
\begin{tabular}{|c|c|ccccc|}
\multicolumn{7}{c}{{\small\bf The Table of HOMFLY polynomials including 4-strand
knots with 6, 7, 8 crossings}}\\
\multicolumn{7}{c}{}\\
\hline
&&&&&&\\
{\rm knot} &$(a_1,b_1,\ldots)$&$S_4^*$&$S_{31}^*$&$S_{22}^*$&
$S_{211}^*$&$S_{1111}^*$\\
&&&&&&\\
\hline
&&&&&&\\
$6_1$&(1,-1,1,1,-1,-2)&$q^{-1}$&$q^{5}-q^{3}-q^{1}+q^{-1}-2q^{-3}+q^{-5}$&$q^{-1}-q^{1}$&
$-q^{5}+2q^{3}-q^{1}+q^{-1}+q^{-3}-q^{-5}$&$(-q)$\\
&&&&&&\\
$7_2$&(-1,1,0,-1,-1,1,0,-1,-3)&$q^{-5}$&$-q^{3}+3q^{1}-4q^{-1}+3q^{-3}-2q^{-5}$&
$-q^{7}+q^{5}-2q^{3}+3q^{1}-3q^{-1}+2q^{-3}-q^{-5}+q^{-7}$&
$2q^{5}-3q^{3}+4q^{1}-3q^{-1}+q^{-3}$&
$-q^5$
\\
&&&&&&\\
$7_4$&(1,-1,0,1,2,-1,0,1,2)&$q^5$&$-3q^{5}+5q^{3}-6q^{1}+5q^{-1}-2q^{-3}$&
$(q-q^{-1})(q^{6}-q^{4}+3q^{2}-1+3q^{-2}-q^{-4}+q^{-6})$&
$2q^{3}-5q^{1}+6q^{-1}-5q^{-3}+3q^{-5}$&
$-q^{-5}$
\\
&&&&&&\\
$7_6$&(-1,1,-1,-1,1,-2)&$q^{-3}$&$-q^{3}+3q^{1}-5q^{-1}+5q^{-3}-4q^{-5}+q^{-7}$&
$-q^{7}+2q^{5}-3q^{3}+3q^{1}-3q^{-1}+3q^{-3}-2q^{-5}+q^{-7}$&
$-q^{7}+4q^{5}-5q^{3}+5q^{1}-3q^{-1}+q^{-3}$&$-q^3$
\\
&&&&&&\\
$7_7$&(1,-1,0,1,-1,1,0,-1,1)&$q$&$-q^{7}+4q^{5}-7q^{3}+7q^{1}-6q^{-1}+2q^{-3}$&
$-q^{7}+3q^{5}-4q^{3}+5q^{1}-5q^{-1}+4q^{-3}-3q^{-5}+q^{-7}$&
$-2q^{3}+6q^{1}-7q^{-1}+7q^{-3}-4q^{-5}+q^{-7}$&$-q^{-1}$
\\
&&&&&&\\
$8_4$&(1,-1,0,1,1,-1,0,1,-3)&$q^{-1}$&$q^{7}-2q^{5}+3q^{3}-4q^{1}+2q^{-1}-q^{-3}-q^{-5}+q^{-7}$&
$(q-q^{-1})(q^4-q^2+1-q^{-2}+q^{-4})$&$-q^{7}+q^{5}+q^{3}-2q^{1}+4q^{-1}-3q^{-3}+2q^{-5}-q^{-7}$&
$-q$
\\
&&&&&&\\
$8_6$&(1,-1,1,1,-1,-4)&$q^{-3}$&$q^{7}-2q^{5}+3q^{3}-3q^{1}+q^{-3}-2q^{-5}+q^{-7}$&$(q-q^{-1})(q^2+1+q^{-2})(q^2-1+q^{-2})$&$-q^{7}+2q^{5}-q^{3}+3q^{-1}-3q^{-3}+2q^{-5}-q^{-7}$&
$-q^3$
\\
&&&&&&\\
$8_{11}$&(1,-1,0,1,-2,1,0,-1,-2)&$q^{-3}$&$(q^2-1+q^{-2})(q^5-q^3-2q^{-3}+q^{-5})$&
$-q^{3}+q^{-3}$&$-q^{7}+3q^{5}-3q^{3}+2q^{1}+q^{-1}-2q^{-3}+2q^{-5}-q^{-7}$&
$-q^3$
\\
&&&&&&\\
$8_{13}$&(1,-1,0,1,2,-1,0,1,-2)&$q$&$-q^{7}+3q^{5}-5q^{3}+4q^{1}-2q^{-1}-q^{-3}+2q^{-5}-q^{-7}$&
$(q-q^{-1})(q^4-q^2+1-q^{-2}+q^{-4})$&$q^{7}-2q^{5}+q^{3}+2q^{1}-4q^{-1}+5q^{-3}-3q^{-5}+q^{-7}$&
$-q^{-1}$
\\
&&&&&&\\
$8_{14}$&(1,-1,0,1,-1,1,0,-1,-3)&$q^{-3}$&$q^{7}-3q^{5}+4q^{3}-3q^{1}+2q^{-3}-3q^{-5}+q^{-7}$&
$(q-q^{-1})(q^2+1+q^{-2})(q^2-1+q^{-2})$&$-q^{7}+3q^{5}-2q^{3}+3q^{-1}-4q^{-3}+3q^{-5}-q^{-7}$&
$-q^3$
\\
&&&&&&\\
$8_{15}$&(1,3,1,1,-1,2)&$q^{-7}$&$q^{3}-4q^{1}+6q^{-1}-7q^{-3}+5q^{-5}-2q^{-7}$&
$(q-q^{-1})(q^{6}-2q^{4}+2q^{2}-3+2q^{-2}-2q^{-4}+q^{-6})$&
$2q^{7}-5q^{5}+7q^{3}-6q^{1}+4q^{-1}-q^{-3}$&$-q^7$
\\
&&&&&&\\
$$\ldots&&&&&&$$\\
\hline
&&&&&&\\
&&$A^3$&$A^1$&$A^{-1}$&$A^{-3}$&\\
&&&&&&\\
\hline
&&&&&&\\
$6_1$&(1,-1,1,1,-1,-2)&1&$-q^{2}+1-q^{-2}$&$-q^{2}+2-q^{-2}$&1&\\
&&&&&&\\
$7_2$&(-1,1,0,-1,-1,1,0,-1,-3)&-1&$q^{2}-1+q^{-2}$&$q^{2}-2+q^{-2}$&$q^{2}-1+q^{-2}$&
\\
&&&&&&\\
$7_4$&(1,-1,0,1,2,-1,0,1,2)&$q^{2}-1+q^{-2}$&$q^{2}-1+q^{-2}$&-1&0&
\\
&&&&&&\\
$7_6$&(-1,1,-1,-1,1,-2)&-1&
$2q^{2}-2+2q^{-2}$&$-q^{4}+2q^{2}-3+2q^{-2}-q^{-4}$&$q^{2}-1+q^{-2}$&
\\
&&&&&&\\
$7_7$&(1,-1,0,1,-1,1,0,-1,1)&
$-q^{2}+2-q^{-2}$&$q^{4}-2q^{2}+4-2q^{-2}+q^{-4}$&$-2q^{2}+2-2q^{-2}$&1&
\\
&&&&&&\\
$8_4$&(1,-1,0,1,1,-1,0,1,-3)&$q^{2}-1+q^{-2}$&$-q^{4}+2q^{2}-2+2q^{-2}-q^{-4}$&
$-q^{4}+q^{2}-2+q^{-2}-q^{-4}$&$q^{2}+q^{-2}$&
\\
&&&&&&\\
$8_6$&(1,-1,1,1,-1,-4)&$q^{2}-1+q^{-2}$&
$-q^{4}+2q^{2}-3+2q^{-2}-q^{-4}$&$-q^{4}+2q^{2}-3+2q^{-2}-q^{-4}$&$q^{2}+q^{-2}$&
\\
&&&&&&\\
$8_{11}$&(1,-1,0,1,-2,1,0,-1,-2)&$q^{2}-1+q^{-2}$&
$-q^{4}+2q^{2}-4+2q^{-2}-q^{-4}$&$-q^{4}+3q^{2}-3+3q^{-2}-q^{-4}$&$q^{2}-1+q^{-2}$&
\\
&&&&&&\\
$8_{13}$&(1,-1,0,1,2,-1,0,1,-2)&$-q^{2}+2-q^{-2}$&
$q^{4}-3q^{2}+4-3q^{-2}+q^{-4}$&$q^{4}-2q^{2}+4-2q^{-2}+q^{-4}$&$-q^{2}+1-q^{-2}$&
\\
&&&&&&\\
$8_{14}$&(1,-1,0,1,-1,1,0,-1,-3)&$q^{2}-2+q^{-2}$&
$-q^{4}+3q^{2}-4+3q^{-2}-q^{-4}$&$-q^{4}+3q^{2}-4+3q^{-2}-q^{-4}$&$q^{2}-1+q^{-2}$&
\\
&&&&&&\\
$8_{15}$&(1,3,1,1,-1,2)&1&$-3q^{2}+2-3q^{-2}$&$2q^{4}-3q^{2}+5-3q^{-2}+2q^{-4}$&
$q^{4}-2q^{2}+3-2q^{-2}+q^{-4}$&
\\
&&&&&&\\
$$\ldots&&&&&&$$\\
\hline
\end{tabular}
}}}

\newpage

\section{5-strand braids}

In the five-strand case one deals with traces over seven representations emerging in the decomposition
\be
[1]^5=[5]+4[41]+5[32]+[311]+[221]+4[2111]+[11111]
\ee
The HOMFLY polynomial has the form
\be
H_{[1]}^{(a_1,b_1,c_1,d_1,\ldots)}
= q^{a_1+b_1+c_1+d_1+\ldots}S_5^* + \\
+ \tr \left[\left(\hat{\mathfrak{R}}_1^{[41]}\right)^{a_1}\ \underbrace{\hat U^{[41]}
\left(\hat{\mathfrak{R}}_1^{[41]}\right)^{b_1}\widetilde{\hat U^{[41]}}}\
\underbrace{\hat U^{[41]}\hat V^{[41]}\hat U^{[41]}
\left(\hat{\mathfrak{R}}_1^{[41]}\right)^{c_1}
\widetilde{\left(\hat U^{[41]}\hat V^{[41]}\hat U^{[41]}\right)}}\times
\right.\\
\left.\times
\underbrace{\hat U^{[41]}\hat V^{[41]}\hat W^{[41]}\hat U^{[41]}\hat V^{[41]}
\hat U^{[41]}
\left(\hat{\mathfrak{R}}_1^{[41]}\right)^{d_1}
\widetilde{\left(\hat U^{[41]}\hat V^{[41]}\hat U^{[41]}\hat W^{[41]}\hat V^{[41]}
\hat U^{[41]}\right)}}\ldots \right]S_{41}^* +\\
+ \tr \left[\left(\hat{\mathfrak{R}}_1^{[32]}\right)^{a_1}\ \underbrace{\hat U^{[32]}
\left(\hat{\mathfrak{R}}_1^{[32]}\right)^{b_1}\widetilde{\hat U^{[32]}}}\
\underbrace{\hat U^{[32]}\hat V^{[32]}\hat U^{[32]}
\left(\hat{\mathfrak{R}}_1^{[32]}\right)^{c_1}
\widetilde{\left(\hat U^{[32]}\hat V^{[32]}\hat U^{[32]}\right)}}\times
\right.\\
\left.\times
\underbrace{\hat U^{[32]}\hat V^{[32]}\hat W^{[32]}
\hat U^{[32]}\hat V^{[32]}\hat U^{[32]}
\left(\hat{\mathfrak{R}}_1^{[32]}\right)^{d_1}
\widetilde{\left(\hat U^{[32]}\hat V^{[32]}\hat U^{[32]}\hat W^{[32]}
\hat V^{[32]}\hat U^{[32]}\right)}}\ldots \right]S_{32}^* + \\
+  \left( q \longrightarrow -\frac{1}{q} \right) \ \ \ \ +\ \ \ \ \alpha S_{311}^*\
\ \ \ \ \ \
\ee
where the contributions of $[5]$ is trivial, those of $[41]$ and $[32]$ can be calculated in the $SU(2)$ case, and
the contributions of $[221]$, $[2111]$ and $[11111]$ are restored by the symmetry. The only non-trivial
contribution that requires the higher group ($SU(3)$) comes from representation $[311]$.

The mixing matrices are described by the following system of the trees:

\vspace{1cm}

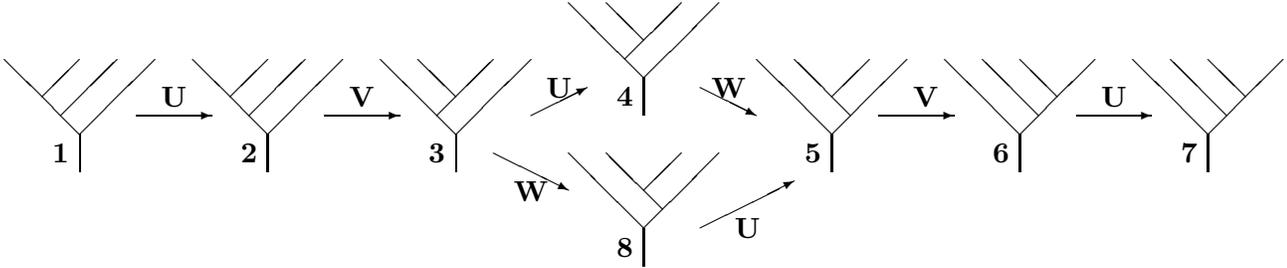
\unitlength=0.25mm
\begin{picture}(20,160)(0,-100)
\put(0,0){\line(0,-1){20}}
\put(0,0){\line(-1,1){40}}
\put(-10,10){\line(1,1){30}}
\put(-20,20){\line(1,1){20}}
\put(-30,30){\line(1,1){10}}
\put(0,0){\line(1,1){40}}
 \put(100,0){\line(0,-1){20}}
\put(100,0){\line(-1,1){40}}
\put(90,10){\line(1,1){30}}
\put(80,20){\line(1,1){20}}
\put(90,30){\line(-1,1){10}}
\put(100,0){\line(1,1){40}}
 \put(200,0){\line(0,-1){20}}
\put(200,0){\line(-1,1){40}}
\put(190,10){\line(1,1){30}}
\put(200,20){\line(-1,1){20}}
\put(190,30){\line(1,1){10}}
\put(200,0){\line(1,1){40}}
 \put(300,30){\line(0,-1){20}}
\put(300,30){\line(-1,1){40}}
\put(290,40){\line(1,1){30}}
\put(300,50){\line(-1,1){20}}
\put(310,60){\line(-1,1){10}}
\put(300,30){\line(1,1){40}}
 \put(400,0){\line(0,-1){20}}
\put(400,0){\line(-1,1){40}}
\put(410,10){\line(-1,1){30}}
\put(400,20){\line(1,1){20}}
\put(410,30){\line(-1,1){10}}
\put(400,0){\line(1,1){40}}
 \put(500,0){\line(0,-1){20}}
\put(500,0){\line(-1,1){40}}
\put(510,10){\line(-1,1){30}}
\put(520,20){\line(-1,1){20}}
\put(510,30){\line(1,1){10}}
\put(500,0){\line(1,1){40}}
 \put(600,0){\line(0,-1){20}}
\put(600,0){\line(-1,1){40}}
\put(610,10){\line(-1,1){30}}
\put(620,20){\line(-1,1){20}}
\put(630,30){\line(-1,1){10}}
\put(600,0){\line(1,1){40}}
%
%
 \put(300,-50){\line(0,-1){20}}
\put(300,-50){\line(-1,1){40}}
\put(310,-40){\line(-1,1){30}}
\put(300,-30){\line(1,1){20}}
\put(290,-20){\line(1,1){10}}
\put(300,-50){\line(1,1){40}}
\put(-10,-10){\makebox(0,0)[cc]{\bf 1}}
\put(90,-10){\makebox(0,0)[cc]{\bf 2}}
\put(190,-10){\makebox(0,0)[cc]{\bf 3}}
\put(290,20){\makebox(0,0)[cc]{\bf 4}}
\put(390,-10){\makebox(0,0)[cc]{\bf 5}}
\put(490,-10){\makebox(0,0)[cc]{\bf 6}}
\put(590,-10){\makebox(0,0)[cc]{\bf 7}}
\put(290,-60){\makebox(0,0)[cc]{\bf 8}}
\put(30,10){\vector(1,0){40}}
\put(50,20){\makebox(0,0)[cc]{\bf U}}
\put(130,10){\vector(1,0){40}}
\put(150,20){\makebox(0,0)[cc]{\bf V}}
\put(240,10){\vector(2,1){30}}
\put(255,25){\makebox(0,0)[cc]{\bf U}}
\put(330,25){\vector(2,-1){30}}
\put(345,25){\makebox(0,0)[cc]{\bf W}}
\put(220,-10){\vector(2,-1){40}}
\put(240,-30){\makebox(0,0)[cc]{\bf W}}
\put(330,-50){\vector(2,1){50}}
\put(355,-50){\makebox(0,0)[cc]{\bf U}}
\put(425,10){\vector(1,0){40}}
\put(450,20){\makebox(0,0)[cc]{\bf V}}
\put(530,10){\vector(1,0){40}}
\put(550,20){\makebox(0,0)[cc]{\bf U}}
\end{picture}

Note that one can pass from tree 3 to tree 5 both through tree 4 and tree 8.
This imposes a consistency condition, which implies that the matrices
$U$ and $W$ should commute.

The $\mathfrak{R}$-matrices in the non-trivial cases are
\be
\hat{\mathfrak{R}}^{[41]} =  \left(\begin{array}{cccc}
q&&&\\&q&&\\&&q&\\&&&-\frac{1}{q}
\end{array}\right), \ \ \ \ \ \ \
\hat{\mathfrak{R}}^{[32]} =  \left(\begin{array}{ccccc}
q&&&&\\&q&&&\\&&-\frac{1}{q}&&\\&&&q&\\&&&&-\frac{1}{q}
\end{array}\right)
\ee
while the matrices $U$, $V$ and $W$ in representation $[41]$ can be read off from (\ref{UV41})-(\ref{W41}) and are
\be
\hat U^{[41]} = \left(\begin{array}{cccc}
1&&&\\
&1&&\\
&&C_2&S_2\\
&&-S_2&C_2
\end{array}\right), \ \ \ \ \
\hat V^{[41]} = \left(\begin{array}{cccc}
1&&&\\
&C_3&S_3&\\
&-S_3&C_3&\\
&&&1
\end{array}\right), \ \ \ \ \
\hat W^{[41]} = \left(\begin{array}{cccc}
C_4&S_4&&\\
-S_4&C_4&&\\
&&1&\\
&&&1
\end{array}\right)
\ee
with
\be
C_k= \displaystyle{{1\over [k]_q}}
\ ,\ \ \ \ \ \ \ \  \
S_k= \displaystyle{{\sqrt{[k]_q^2-1}\over [k]_q}}
\ee

Similarly,
\be
\hat U^{[32]} = \left(\begin{array}{ccccc}
1&&&&\\
&-C_2&S_2 &&\\
&-S_2&-C_2&&\\
&&&-C_2&S_2\\
&&&-S_2&-C_2
\end{array}\right), \ \ \ \ \ \
\hat V^{[32]} = \left(\begin{array}{ccccc}
-C_3&S_3&&&\\
S_3&C_3& &&\\
&&1&&\\
&&&1&\\
&&&&-1
\end{array}\right), \\
\hat W^{[32]} = \left(\begin{array}{ccccc}
1&&&&\\
&-C_2&&S_2 &\\
&&-C_2&&S_2\\
&S_2&&C_2&\\
&&S_2&&C_2
\end{array}\right)
\ee

From these manifest matrices one can check that, indeed,
\be
\l[\hat U^{[32]},\hat W^{[32]}] = 0\hspace{1cm}\hbox{and}\hspace{1cm}
\l[\hat U^{[41]},\hat W^{[41]}] = 0
\ee

\newpage

\thispagestyle{empty}

\hspace{-2.cm}
{\tiny{\rotate{
\begin{tabular}{|c|cccccc|}
\multicolumn{7}{c}{{\small\bf The Table of HOMFLY polynomials including three 5-strand knots with 8 crossings,
$8_1,\ 8_3,\ 8_{12}$ and first 3 torus knots with 5 strands}}\\
\multicolumn{7}{c}{}\\
\hline
&&&&&&\\
{\rm knot} & $(a_1,b_1,c_1,d_1|a_2,b_2,c_2,d_2|a_3,b_3,c_3,d_3,\ldots)$ &
$A^4$ & $A^2$ & 1 & $A^{-2}$ & $A^{-4}$ \\
&&&&&&\\ \hline
&&&&&&\\
$8_1$&$(-2,-1,0,0|1,-1,-1,1|0,1,-1,1)$:&
1&$-q^2+1-q^{-2}$&$-(q-q^{-1})^2$&$-(q-q^{-1})^2$&1\\
&&&&&&\\
$8_3$& $(-2,-1,1,0|1,-1,1,1|0,0,-1,1)$:&
1&$-(q-q^{-1})^2$&$-2q^2+3-2q^{-2}$&$-(q-q^{-1})^2$&1\\
&&&&&&\\
$8_{12}$ & $(-1,1,-1,1|-1,1,-1,1)$:&
1&$-2q^2+3-2q^{-2}$&$q^{4}-3q^{2}+5-3q^{-2}+q^{-4}$&$-2q^2+3-2q^{-2}$&1\\
&&&&&&\\ \ldots
&&&&&&\\
&&&&&&\\
\hline
&&&&&&\\
$\l[5,1]$ \ {\rm unknot}  &$(1,1,1,1)$:&1&&&&\\
&&&&&&\\
$\l[5,2]=[2,5]$  &$(-1,-1,-1,-1|-1,-1,-1,-1)$:&&&&$-[4]_q/[2]_q$&
$[6]_q/[2]_q$\\
&&&&&&\\
$\l[5,3]=[3,5]$ &$(-1,-1,-1,-1|-1,-1,-1,-1|-1,-1,-1,-1)$:&&&$[4]_q/[2]_q$&
$[4]_q[2]_q(q^2-1+q^{-2})$&$[7]_q(q^5-q^4+q-1+q^{-1}-q^{-4}+q^{-5})$\\
&&&&&&\\ \ldots
&&&&&&\\
\hline
&&&&&&\\
{\rm knot} &$(a_1,b_1,c_1,d_1|a_2,b_2,c_2,d_2|a_3,b_3,c_3,d_3,\ldots)$ &
$S_5^*$& \multicolumn{2}{c}{$S_{41}^*$} & \multicolumn{2}{c|}{$S_{32}^*$}\\
&&&&&&\\ \hline
&&&&&&\\
$8_1$&$(-2,-1,0,0|1,-1,-1,1|0,1,-1,1)$:&
$q^{-2}$&\multicolumn{2}{c}{$2q^{4}-4q^{2}+2-2q^{-4}+q^{-6}$}&
\multicolumn{2}{c|}{$(q-q^{-1})\frac{q^{10}-2q^8+3q^6-4q^4+2q^2-1}{q^3}$}
\\
&&&&&&\\
$8_3$&$(-2,-1,1,0|1,-1,1,1|0,0,-1,1)$:&
1&\multicolumn{2}{c}{$q^{6}-q^{4}-q^{2}+1-q^{-2}-q^{-4}+q^{-6}$}&
\multicolumn{2}{c|}{$(q^2-1+q^{-2})(q-q^{-1})^2$}
\\
&&&&&&\\
$8_{12}$ &$(-1,1,-1,1|-1,1,-1,1)$:&
1&\multicolumn{2}{c}{$2q^{4}-6q^{2}+7-6q^{-2}+2q^{-4}$}&
\multicolumn{2}{c|}{$[4]_q/[2]_q(q-q^{-1})^2(q^2-1+q^{-2})^2$}
\\
&&&&&&\\ \ldots
&&&&&&\\
&&&&&&\\
\hline
&&&&&&\\
$\l[5,1]$ \ {\rm unknot}  &$(1,1,1,1)$:&$q^4$&\multicolumn{2}{c}{$-q^2$}&
\multicolumn{2}{c|}{0}\\
&&&&&&\\
$\l[5,2]=[2,5]$  &$(-1,-1,-1,-1|-1,-1,-1,-1|-1,-1,-1,-1)$:&$q^{-8}$&
\multicolumn{2}{c}{$-q^{-4}$}&\multicolumn{2}{c|}{0}\\
&&&&&&\\
$\l[5,3]=[3,5]$ &$(-1,-1,-1,-1|-1,-1,-1,-1|-1,-1,-1,-1)$:&$q^{-12}$&
\multicolumn{2}{c}{$q^{-6}$}&\multicolumn{2}{c|}{0}\\
&&&&&&\\ \ldots
&&&&&&\\
\hline
&&&&&&\\
{\rm knot} &
\multicolumn{2}{c}{$S_{311}^*$} & \multicolumn{2}{c}{$S_{221}^*$} & $S_{2111}^*$ & $S_{11111}^*$\\
&&&&&&\\ \hline
&&&&&&\\
$8_1$&
\multicolumn{2}{c}{$-2q^{6}+4q^{4}-4q^{2}+5-4q^{-2}+4q^{-4}-2q^{-6}$}&
\multicolumn{2}{c}{$(q-q^{-1})(q^{3}-2q+4q^{-1}-3q^{-3}+2q^{-5}-q^{-7})$}&
$q^{6}-2q^{4}+2-4q^{-2}+2q^{-4}$&$q^2$
\\
&&&&&&\\
$8_3$&
\multicolumn{2}{c}{$-2q^{6}+3q^{4}-q^{2}+1-q^{-2}+3q^{-4}-2q^{-6}$}
&\multicolumn{2}{c}{$(q^2-1+q^{-2})(q-q^{-1})^2$}&
$q^{6}-q^{4}-q^{2}+1-q^{-2}-q^{-4}+q^{-6}$&1
\\
&&&&&&\\
$8_{12}$ &
\multicolumn{2}{c}{$q^{8}-6q^{6}+11q^{4}-12q^2+13-12q^{-2}+11q^{-4}
-6q^{-6}+q^{-8}$}
&\multicolumn{2}{c}{$[4]_q/[2]_q(q-q^{-1})^2(q^2-1+q^{-2})^2$}&
$2q^{4}-6q^{2}+7-6q^{-2}+2q^{-4}$&1
\\
&&&&&&\\ \ldots
&&&&&&\\
&&&&&&\\
\hline
&&&&&&\\
$\l[5,1]$ \ {\rm unknot}  &\multicolumn{2}{c}{1}&\multicolumn{2}{c}{0}&$-q^{-2}$&$q^{-4}$\\
&&&&&&\\
$\l[5,2]=[2,5]$  &\multicolumn{2}{c}{1}&\multicolumn{2}{c}{0}&$-q^4$&$q^8$\\
&&&&&&\\
$\l[5,3]=[3,5]$ &\multicolumn{2}{c}{1}&\multicolumn{2}{c}{0}&$q^6$&$q^{12}$\\
&&&&&&\\ \ldots
&&&&&&\\
\hline
\end{tabular}
}}}

\section{Summary and comments}

The goal of this paper was to present an expansion of the HOMFLY polynomials of knots into
characters (the Schur functions)
\be
H_R^{\cal K} = \sum_{Q_{\vdash} m|R|} h_R^Q S_Q^* \label{Hexpc}
\ee
{\it A priori}, such an expansion is ambiguous for $m>3$, and we
suggest to define it through the braid realization ${\cal B}$ of
the knot ${\cal K}$, which contains $m$ strands, and present
explicit formulas for $h_R^Q$ for $R=[1]$ and $m=2,3,4$. In the case
of $m=5$ one of the coefficients, $h_{[1]}^{[311]}$ is left
theoretically undetermined in the present paper, and we read off its values
for $5$-strand knots from the known answers for the HOFMLY polynomials.

The answer for a fixed braiding representation of the knot
is completely expressed through the known eigenvalues of the quantum
${\cal R}$-matrices and through the mixing ${\cal U}$-matrices
 that relates ${\cal R}$-matrices acting on
different strands of the braid. The matrix elements of these mixing matrices
are the Racah coefficients, they turn out to look very simple and
exhibit vast universality properties.

More concretely, for the $m$-strand braid with the representation $R$ running through
each strand, we worked in the basis of irreps at each transverse slice of the braid,
$R^{\otimes m}\to \oplus Q$. There were a few crucial ingredients in the construction.
\begin{itemize}
\item Typically, there emerge several similar irreps in this expansion.
They are labeled by the path leading to the concrete irrep through the fusion tree $T$.
Since we considered here only the trees with the vertices where at least one of the
edges carries $R=[1]$, we can describe the path through the tree
only with the sequence of non-trivial entries: $V_T(R_1\to R_2\to R_3...)$.
\item Thus, typically there is a set of similar irreps for the given tree. These sets
for different trees are connected by the mixing matrices. We require the
vectors in the irrep to be normalized. Then, the mixing matrices are orthogonal.
\item If there are $N$ different similar irreps, the mixing matrix is of the size $N\times
N$. However, we demonstrated that in the considered cases it is really $2\times 2$ matrix,
because of the specific structure of mixing. It is illustrated in the Table, where
the representations mixing are marked by the same letter:

\begin{tabular}{ccccc}
{\rm paths} & & U & V & W \\
\hline
&{\rm level} & 1 & 2 & 3 \\
&{\rm mixing\ reps} & [2]\&[11] & [3]\&[21] & [4]\&[31] \\
&&&& [31]\&[22] \\
\hline
$V_T([2]\to [3]\to [4]\to [41])$ &&&&z\\
$V_T([2]\to [3]\to [31]\to [41])$ &&&y&z\\
$V_T([2]\to [21]\to [31]\to [41])$ &&x&y&\\
$V_T([11]\to [21]\to [31]\to [41])$ &&x&&\\
\hline
$V_T([2]\to [3]\to [31]\to [32])$ &&&c&\\
$V_T([2]\to [21]\to [31]\to [32])$  &&a&c&e\\
$V_T([11]\to [21]\to [31]\to [32])$ &&a&&d\\
$V_T([2]\to [21]\to [22]\to [32])$ &&b&&e\\
$V_T([11]\to [21]\to [22]\to [32])$ &&b&&d \\
\hline
\end{tabular}

\bigskip

One can see that only the paths that differ from each other exactly in one item can mix, which
effectively leads to $2\times 2$ matrices.
\item Actually, the mixing matrices ${\cal U}_{\mu\nu}$
have their own hierarchical structure,
expressed as a relation between two decorated routed trees,
decorated by the
representation in the decomposition
\be
R^{\otimes m} = \sum_Q {\cal M}_Q\otimes Q
\ee
${\cal M}_Q$ is actually labeled by a 3-valent tree,
and ${\cal U}$ matrix is decomposed into a product of elementary
constituents $(U,V,W,\ldots)$, realizing elementary steps, that is,
flips of a single edge in the routed 3-valent tree.
\item In order to calculate the knot invariant one also needs to know the eigenvalues of
${\cal R}$-matrix. Then, one chooses the ${\cal R}$-matrix acting at the first two strands
diagonal and expresses ${\cal R}$-matrices acting at other strands through diagonal
${\cal R}$-matrix rotated with the mixing matrices. It immediately gives the knot polynomial.
\item Unfortunately, the coefficients of expansion (\ref{Hexpc}) $h_R^Q$ are not
knot invariant. For instance, the HOMFLY polynomial for toric knot $T[n,m]$ is the same as
that for toric knot $T[m,n]$, but the numbers of strands and, therefore, expansions (\ref{Hexpc})
are different for them. However, it looks plausible that while the full knot-{\it in}variance
is lost in transfer to $h_R^Q$, some rich knot-{\it co}variance w.r.t. switches between
different braid realizations can be finally found for the character expansion (\ref{Hexpc}).
\end{itemize}
Among the advantages of the character expansion (\ref{Hexpc})
we emphasize the following:
\begin{itemize}
\item Its coefficients do not depend
on $A$, which allows one to calculate them for small rank groups $SU_q(l)$, with
needed $l$ depending only on $R$ and $m$, still knot invariants are obtained at once for
all the values of the rank of the {\it gauge} group $SU(N)$.
Moreover, the answers appear directly expressed through $A=q^N$
as required for the HOMFLY polynomials (what is not at all automatic
for other calculations in Chern-Simons theory). This gives a good technical alternative to the
usually explored in calculating skein relations.
\item There are also two more advantages of more theoretical nature (see \cite{AMMceI}).
First of all, the character expansion can be naturally extended to
new set of variables: one can consider the Schur functions at
arbitrary points. This opens a way to dealing with new set of
questions related, e.g., to integrability. Second, one can
effectively use this expansion in order to immediately construct
(quantum) $A$-polynomials.
\end{itemize}

We plan to develop this formalism
and its applications
in the forthcoming publications.

\section*{Acknowledgements}

Our work is partly supported by Ministry of Education and Science of
the Russian Federation under contract 14.740.11.081, by RFBR
grants 10-02-00509 (A.Mir.), 10-02-00499 (A.Mor.) and 10-02-01315 (And.Mor.),
by joint grants 11-02-90453-Ukr, 09-02-93105-CNRSL, 09-02-91005-ANF,
10-02-92109-Yaf-a, 11-01-92612-Royal Society.


\begin{thebibliography}{12}

\bibitem{charexp} D.J.Gross and W.Taylor, Nucl.Phys. {\bf B400} (1993) 181-210, arXiv:hep-th/9301068

\bibitem{cft} J.Cardy, Nucl.Phys. {\bf B270} (1986) 186-204\\
E.Melzer, Lett.Math.Phys. {\bf 31} (1994) 233-246, hep-th/9312043

\bibitem{mamocharexp} C.Itzykson and J.Zuber, J.Math.Phys. {\bf 21} (1980) 411\\
V.Kazakov, M.Staudacher and T.Wynter,
Commun.Math.Phys. {\bf 177} (1996) 451-468

\bibitem{intcharexp} M. Jimbo, T. Miwa, Publ.RIMS, Kyoto Univ,. {\bf 19} (1983) 943-1001\\
S.Kharchev, A.Marshakov, A.Mironov and A.Morozov, Int. J. Mod. Phys. {\bf A10} (1995) 2015,
hep-th/9312210\\
A.Alexandrov, A.Mironov, A.Morozov and S.Natanzon, arXiv:1103.4100

\bibitem{latYM}
E.Brezin and D.Gross, 
Phys.Lett., {\bf B97} (1980) 120\\
D.Gross and E.Witten, 
Phys.Rev., {\bf D21} (1980) 446-453\\
B.De Wit and G.t'Hooft, 
Phys.Lett., {\bf B69} (1977) 61

\bibitem{mamoYM}
A.Mironov, A.Morozov and G.Semenoff, 
Int.J.Mod.Phys., {\bf A10} (1995) 2015, hep-th/9404005\\
A.Morozov,
Theor.Math.Phys. {\bf 162} (2010) 1-33 (Teor.Mat.Fiz. {\bf 161}
(2010) 3-40),
arXiv:0906.3518\\
A.Balantekin, arXiv:1011.3859\\
A.Mironov, A.Morozov and Sh.Shakirov,
JHEP {\bf 1103} (2011) 102, arXiv:1011.3481

\bibitem{AMMceI} A.Mironov, A.Morozov and And.Morozov, arXiv:1112.5754

\bibitem{CS} E.Witten,
Commun. Math. Phys. {\bf 121} (1989) 351-399\\
R.H.Kaul, Commun.Math.Phys. {\bf 162} (1994) 289-320, hep-th/9305032\\
Zodinmawia and P.Ramadevi, arXiv:1107.3918

\bibitem{HOMFLY} P.Freyd, D.Yetter, J.Hoste, W.B.R.Lickorish, K.Millet, A.Ocneanu,
Bull. AMS. {\bf 12} (1985) 239\\
J.H.Przytycki and K.P.Traczyk, 
Kobe J. Math. {\bf 4} (1987) 115-139

\bibitem{GSV} S.Gukov, A.Schwarz and C.Vafa, Lett.Math.Phys. {\bf 74} (2005) 53-74,
arXiv:hep-th/0412243

\bibitem{DMMSS}
P.Dunin-Barkowski, A.Mironov, A.Morozov, A.Sleptsov and A.Smirnov, arXiv:1106.4305

\bibitem{knotcharexp}
M.Rosso and V.F.R.Jones, J. Knot Theory Ramifications, {\bf 2} (1993) 97-112\\
J.M.F.Labastida and M.Marino,
J.Knot Theory Ramifications, {\bf 11} (2002) 173\\
X.-S.Lin and H.Zheng, Trans. Amer. Math. Soc. 362 (2010) 1-18 math/0601267

\bibitem{Apoly}
R.Gelca, Math. Proc. Cambridge Philos. Soc. {\bf 133} (2002) 311-323\\
R.Gelca and J.Sain, J. Knot Theory Ramifications, {\bf 12} (2003) 187-201\\
S.Garoufalidis and T.Le, Geometry and Topology, {\bf 9} (2005) 1253-1293, math/0309214

\bibitem{TR}
E.Guadagnini, M.Martellini and M.Mintchev, In Clausthal 1989, Proceedings, Quantum groups", 307-317;
Phys.Lett. B235 (1990) 275\\
N.Yu.Reshetikhin and V.G.Turaev, 
Comm. Math. Phys. {\bf 127} (1990) 1-26

\bibitem{MorSmi}
A.Morozov and A.Smirnov,
Nucl.Phys.B 835:284-313, 2010, arXiv:1001.2003\\
A.Smirnov, hep-th/0910.5011

\bibitem{MF} B.Bakalov and A.Kirillov Jr., {\it Lectures on tensor categories and modular functors},
University Lecture Series 21, American Mathematical Society, Providence, RI
2001\\
V.Turaev and O.Viro, 
Topology, {\bf 31} (1992) 865-902\\
J.W.Barrett and B.W.Westbury, 
Trans. Amer. Math. Soc. {\bf 348} (1996) 3997-4022

\bibitem{katlas}
Knot Atlas at http://katlas.org/wiki/Main\_Page (by D.Bar-Natan)

\bibitem{hook} William Fulton, {\it Young Tableaux, with Applications to
Representation Theory and Geometry}, Cambridge University Press, 1997

\bibitem{CJ} A.Mironov, A.Morozov and S.Natanzon, Theor.Math.Phys. 166 (2011) 1-22,
arXiv:0904.4227; Journal of Geometry and Physics, {\bf 62} (2012) 148-155
arXiv:1012.0433

\bibitem{torus}
S.Stevan, Annales Henri Poincare, {\bf 11} (2010) 1201-1224, arXiv: 1003.2861\\
A.Brini, B.Eynard and M.Marino, arXiv:1105.2012

\bibitem{KhR} M.Khovanov and L.Rozhansky, Fund. Math. {\bf 199} (2008) 1, math.QA/0401268;
Geom. Topol. {\bf 12} (2008) 1387, math.QA/0505056

\bibitem{DGR}
N.M.Dunfield, S.Gukov and J.Rasmussen, Experimental Math. 15 (2006) 129-159, math/0505662\\
E.Gorsky, arXiv:1003.0916\\
M.Aganagic and Sh.Shakirov, arXiv: 1105.5117\\
N.Carqueville and D.Murfet, arXiv:1108.1081



\end{thebibliography}
\end{document}